\documentclass[10pt]{article}
\usepackage{geometry}
\usepackage{cancel,amssymb,amsmath,mathrsfs}
\usepackage{color,tikz}
\usetikzlibrary{patterns, decorations.markings,arrows, calc}
\usepackage{framed}
\usepackage[colorlinks=true, pdfstartview=FitV, linkcolor=darkblue, citecolor=darkblue, urlcolor=darkblue]{hyperref}
\definecolor{shadecolor}{rgb}{0.98, 0.98, 0.9}
\definecolor{darkgreen}{rgb}{0.2, 0.5,  0}
\definecolor{darkblue}{rgb}{0.1,0.1,0.45}
\definecolor{red}{rgb}{0.9,0,0}

\def \bea#1\eea {\begin{align} #1 \end{align}}

\def\&{\hspace{-15pt}&}
\def\nn{\nonumber}
\def \scr{\mathscr}

\def \h{\mathfrak h}

\def\ds{\displaystyle}

\def\ov{\overline}

\def \eqref#1{(\ref{#1})}

\def \remove #1 { \sout{ {\color{red} #1}}}
\def\Re{\mathrm {Re}\,}
\def\Im{\mathrm {Im}\,}
\def \wt{\widetilde}
\def\Cauchy{ \mathbf {C}}
\def \wh{\widehat}
\newcommand{\C}{\mathbb{C}}
\renewcommand{\le}{\left}
\newcommand{\ri}{\right}
\newcommand{\R}{\mathbb{R}}
\newcommand{\Z}{\mathbb{Z}}

\newcommand{\1}{\mathbf{1}}

\renewcommand{\d}{\mathrm d}
\renewcommand{\mod}{\,\mathrm{mod}\,}
\newcommand{\pa}{\partial}

\def\res{\mathop{\mathrm {res}}\limits_}
\def \tr {\mathrm{tr}\,}
\def\be{\begin{equation}}
\def\ee{\end{equation}}

\def\bg{\begin{gathered}}
\def\eg{\end{gathered}}

\newtheorem{theorem}{Theorem}[section]

\newtheorem{lemma}[theorem]{Lemma}
\newtheorem{remark}[theorem]{Remark}
\newtheorem{problem}[theorem]{Riemann--Hilbert Problem}

\newtheorem{proposition}[theorem]{Proposition} 
\newtheorem{corollary}[theorem]{Corollary} 
 
\newtheorem{definition}[theorem]{Definition}
\newtheorem{assumption}[theorem]{Assumption}
\def\bet
{
\begin{shaded}
\begin{theorem}
}
\def\eet{\end{theorem} \end{shaded}}

\def\bp
{
\begin{shaded}
\begin{proposition}
}
\def\ep{
\end{proposition}\end{shaded}
}

\def\QED {\hfill $\blacksquare$\par\vskip 10pt}
\def\s{\sigma}

\def\E{\mathcal{E}}
\def \DD{\mathbb D}

\renewcommand{\theequation}{\arabic{section}.\arabic{equation}}

\makeatletter
\@addtoreset{equation}{section}
\makeatother

\def\div{\mathrm{div}\,}

\begin{document}
\vspace{0.2cm}
\begin{center}
\begin{Large}
\textbf{Nonlinear steepest descent approach to orthogonality on elliptic curves} 
\end{Large}
\end{center}

\begin{center}
M. Bertola$^{\dagger\ddagger\diamondsuit}$ \footnote{Marco.Bertola@\{concordia.ca, sissa.it\}}, 
\\
\bigskip
\begin{minipage}{0.7\textwidth}
\begin{small}
\begin{enumerate}
\item [${\dagger}$] {\it  Department of Mathematics and
Statistics, Concordia University\\ 1455 de Maisonneuve W., Montr\'eal, Qu\'ebec,
Canada H3G 1M8} 
\item[${\ddagger}$] {\it SISSA, International School for Advanced Studies, via Bonomea 265, Trieste, Italy }
\item[${\diamondsuit}$] {\it Centre de recherches math\'ematiques,
Universit\'e de Montr\'eal\\ C.~P.~6128, succ. centre ville, Montr\'eal,
Qu\'ebec, Canada H3C 3J7}
\end{enumerate}
\end{small}
\end{minipage}
\vspace{0.5cm}
\end{center}

\begin{abstract}
We consider the recently introduced notion of denominators of Pad\'e--like approximation problems on a Riemann surface. These denominators are  related as in the classical case to the notion of orthogonality over a contour.  We investigate a specific setup where the Riemann surface is a real elliptic curve with two components and the measure of orthogonality is supported on one of the two real ovals. Using a characterization in terms of a Riemann--Hilbert problem, we determine the strong asymptotic behaviour of the corresponding orthogonal functions for large degree. The theory of vector bundles and the non-abelian Cauchy kernel play a prominent role even in this simplified setting, indicating the new challenges that the steepest descent method  on a Riemann surface has to overcome. 
\end{abstract}

\setcounter{tocdepth}{12}
\tableofcontents

\section{Introduction}
Orthogonal polynomials   have a plethora of applications in modern mathematics ranging from the pure approximation theory aspect (linked to their origins as denominators in the Pad\'e\ approximation problem) \cite{Baker, Szegobook}, signal analysis \cite{Grunbaum1}, quantum mechanics, random matrices \cite{Deift}, combinatorics \cite{Stanton}.

Orthogonal polynomials  are functions on the Riemann sphere $\mathbb P^1$ with a single pole at the point $\infty$: we can thus interpret them as sections of the line bundles $\mathcal O(n\infty)$, for $n=0,1,\dots$. 

In a recent work \cite{BertoPade} this point of view was exploited to generalize the notion of orthogonal polynomials to the setting where, instead of the Riemann sphere $\mathbb P^1$, we consider meromorphic functions on an arbitrary Riemann surface $\mathcal C$ of genus $g\geq 1$. 
The new notion introduced in loc. cit. bears the same relationship to the Pad\'e\ approximation problem as the usual notion of  orthogonal polynomials, with the appropriate understanding of the Weyl-Stieltjes function as a differential  on $\mathcal C$ with prescribed singularities.

An important problem, both in and by itself as well as for the applications, is that of describing the asymptotic behaviour of the orthogonal polynomials when the degree tends to infinity \cite{DKMVZ, BleherIts}. 

The approaches to the study of asymptotics use the integral representation of classical orthogonal polynomials,  or tools from potential theory that allow to cast a wider net but with weaker results in general \cite{SaffTotik}, or  tools of integrable systems. These latter  provide complete control over  strong asymptotic results and  have been honed over the past two decades or so since their initial use in \cite{DKMVZ, BleherIts}. They are referred to interchangeably  as "Riemann--Hilbert method" or "Deift-Zhou method" or, as we do in the title, "nonlinear steepest descent analysis".

The main goal of the present paper is to apply those techniques to the new notion of orthogonality on a Riemann surface. A similar setup was recently considered in \cite{olver} but with some different definitions of  the basis of functions that are orthogonalized:
{  
in that case the authors consider a (real) cubic curve $y^2 = P_3(x)$ (with $P_3$ a polynomial of degree $3$) and a measure on its real component (i.e. as a subset of $\R^2$). Then they investigate the definition and properties of the orthogonalization of the (ordered) set of  ``monomials'' $\{1, x, x^2, y, xy, x^2y, \dots\}$ where the pairing is the integration on the real cubic curve with respect to the given measure. In our present paper we apply the orthogonalization procedure to a  different set of ``monomial'' functions which are related to a Pad\'e\ problem as introduced in \cite{BertoPade} and also expressible as a suitable Riemann--Hilbert problem (described in detail below). 
}

The present result is, to the knowledge of the author, the first example of nonlinear steepest descent analysis for orthogonality on a Riemann surface,  and new obstacles and corresponding solutions tie the method inextricably to the theory of vector bundles, as we discuss later. 

Specifically we aim at describing the asymptotic behaviour of the orthogonal ``polynomials'', for large degree, as functions on the Riemann surface. This includes their oscillatory behaviour in the neighbourhood of the support of the orthogonality  measure, their behaviour away from the support, the asymptotic density of their zeros. 

The path  was opened in \cite{BertoPade} by setting up a Riemann--Hilbert problem  on the Riemann surface $\mathcal C$ which characterizes the orthogonal ``polynomials'' in terms of the solution of a matrix--valued function. This is the analog of the standard formulation of orthogonal polynomials by Fokas-Its-Kitaev \cite{FIK} that paved the way to the analysis of Bleher and Its in \cite{BleherIts} (for particular class of measures of orthogonality) and Deift-Kriecherbauer-McLaughlin-Venakides-Zhou in  \cite{DKMVZ} (for a wider class).\par\vskip 3pt

In \cite{BertoPade} it was not clear whether the new notion of these orthogonal ``polynomials'' would be amenable to the steepest descent analysis and hence whether the Riemann--Hilbert formulation  had a practical usefulness beyond the mere elegance of the formulation.  

The goal of the present paper is to show a class of problems where the above doubt is resolved favourably; namely we consider the orthogonality of ``polynomials'' on a real  elliptic curve $\mathcal C= \mathcal E_\tau$ (i.e. a Riemann surface of genus $1$ with an anti holomorphic involution). The orthogonality depends on  an arbitrary positive measure on a connected component   of the real locus.

As customary in these sorts of analyses, there are two possible setups; either the measure of orthogonality is fixed, or it depends on $n$ (``scaling regime'') in the sense that the density with respect to a fixed measure (e.g. the Lebesgue one on the real axis in genus zero) is raised to a power that is scaled together with the degree of the polynomial.

In this paper we consider the first case, where the measure is fixed; we will comment in Section \ref{final} on how the present results indicate the steps for an analysis also in the second case, and what is needed to proceed in that case.

\subsection{Description of setup and results}
We first summarize briefly the notion and notations from \cite{BertoPade} and then progressively specialize.\\[1pt]
 Let $\mathcal C$ be a smooth algebraic curve of genus $g\geq 1$ and $\gamma\subset \mathcal C$ a smooth contour (or collection thereof) with a smooth measure $\d\mu$ on it: by this we mean that it can be written as $f(z)\d z$ in a local coordinate with $f(z)$ a smooth function (at this stage it could be complex valued).
 Let $\infty\in \mathcal C\setminus \gamma$ be a chosen point, and $\scr D$ a nonspecial divisor in $\mathcal C\setminus (\gamma \cup\{\infty\})$. Namely,  $\scr D$ consists of $g$ points (counted with multiplicity) and such that there is no non-constant meromorphic function on $\mathcal C$ with poles only at points of $\scr D$ (of the order bounded by the  corresponding multiplicity of the point). We recall (Riemann--Roch theorem) that this is the {\it generic} situation. For genus $g=1$ there are no special divisors of degree $1$ and hence the issue is moot.

We now need to describe what is the replacement of the notion of ``orthogonal polynomials''. Instead of polynomials we consider appropriate meromorphic functions, characterized by suitable requirements  on the position and degree of their poles, thus mimicking the fact that a polynomial of degree $n$ is simply a meromorphic function with a pole of order $n$ at infinity. 

In the algebro--geometric setup such meromorphic functions are identified with {\it sections} of line bundles and therefore we will use the terminology {\it orthogonal sections}, while hoping that this terminology will not be too distracting for the reader not attuned to these notions.
Thus we consider the spaces
\be
\scr P_n:= \mathcal O(n\infty + \scr D) = \le\{f:\mathcal C\to \mathbb P^1:\ \ \ \div(f)\geq - n\infty-\scr D \ri\}.
\ee
These are the analogs of the space of polynomials of degree $n$: in particular $\dim_\C \scr P_n = n+1$, as  follows from the Riemann--Roch theorem.

The orthogonal sections are thus defined 
\begin{definition}
The orthogonal sections are a sequence $\{\pi_n\}_{n\in \mathbb N}$ with $\pi_n\in \scr P_n\setminus \scr P_{n-1}$  such that
\be
\int_\gamma \pi_n(p) \pi_m(p) \d\mu(p) = h_n \delta_{nm}.
\ee
\end{definition}
It is evident that even if such a sequence exists, we can produce another one  by a rescaling $\pi_n(p)\to c_n\pi_n(p)$, with $c_n\neq 0$ an arbitrary sequence of nonzero numbers. 
To dispose of this arbitrariness we assume further that we have chosen a local coordinate $z$ near $\infty$ (such that $z(\infty)=0$) and we fix their normalization by the requirement 
\be
\pi_n(p) = \frac 1{z^n(p)}  \le(1+ \mathcal O(z(p))\ri).
\ee
We refer to these as {\it monic} orthogonal sections (relative to the choice of coordinate). 

It was shown in \cite{BertoPade} that these orthogonal sections can be characterized as the $(1,1)$ entry of the solution of the following Riemann--Hilbert-Problem (RHP)
\begin{problem}
\label{introYRHP}
Let $Y(p)$ be a $2\times 2$ matrix defined on $\mathcal C\setminus \gamma$ such that the first column consists of meromorphic functions and the second column consists of meromorphic differentials satisfying the following properties:
\begin{enumerate}
\item $Y(p)$ admits bounded boundary values that satisfy 
\be
Y(p_+) = Y(p_-) \le[\begin{array}{cc}
1 & \d\mu(p)\\
0 & 1
\end{array}\ri],\ \ \ \ \ p\in \gamma.
\ee
\item The entries have the  divisor properties (away from $\gamma$) as 
\be
Y(p) = \le[\begin{array}{cc}
\mathcal O\le(n\infty+\scr D\ri) & \mathcal K\le(-(n-1)\infty - \scr D\ri)\\
\mathcal O\le((n-1)\infty+\scr D\ri) &  \mathcal K\le(-(n-2)\infty - \scr D\ri)
\end{array}\ri]
\ee
\item In the local coordinate $z=z(p)$ we have the normalization  near $\infty$ (so that $z(\infty)=0$) 
\be
Y_{11}(p) = \frac 1{z^n} + \mathcal O(z^{-n+1}), \ \ \ Y_{22}(p) = z^n\frac {\d z}{z^{2}}\big(1 + \mathcal O(z) \big)
\ee
\end{enumerate}
\end{problem}
(We point at the paragraph "Notations" at the end of the introduction for an explanation of the symbols $\mathcal O, \mathcal K$)
{  
 We are using a slightly different normalization from the one used in \cite{BertoPade}: if $\wt Y$ is the solution to the RHP 2.8 of loc. cit, then 
$$
Y(p) = \le[\begin{array}{cc}
-i&0\\
0 & -2\pi
\end{array}\ri]\wt Y(p)\le[\begin{array}{cc}
i &0\\
0 &  \frac {1}{2\pi}
\end{array}\ri].
$$}
It was also shown in loc.cit. that the solution to the RHP \ref{introYRHP} exists and is unique if and only if the following determinant does not vanish:
\be
D_n:= \det\big[\mu_{ab}\big]_{a,b=0}^{n-1},\ \ \ \ \mu_{ab} = \int_{\gamma} \sigma_a(p) \sigma_b(p)\d\mu(p),
\ee
where $\sigma_j(p) \in \scr P_j\setminus \scr P_{j-1}$, $j=0,1,\dots$, is an arbitrary sequence of sections. 

{  
 Let us briefly comment on the choice of $\scr D$ and why we impose that  the point $\infty$ does not belong to the divisor. The main reason is in the origin of the problem as Pad\'e\ approximation; indeed (see \cite{BertoPade} for more details)  the orthogonal sections (meromorphic functions)  $\pi_n$ are also interpretable as denominators of the approximation of the ``Weyl-Stieltjes differential"
\be
\mathcal W(p) = \int_{q\in\gamma} C(p, q) \d\mu(q)
\ee
where here $C(p,q)$ is the (unique) Cauchy kernel with the properties  that 
\begin{enumerate}
\item it is a differential in the variable $p$ and a meromorphic function in the variable $q$;
\item $\div_p(C(p,q)) \geq -\infty+\scr D-q$
\item $\div_q(C(p,q)) \geq +\infty -\scr D -p$;
\item $\res{p=q} C(p,q)=1 = -\res{p=\infty} C(p,q)$.
\end{enumerate}
If $\infty \in \scr D$ then, effectively, we have a positive divisor $\wt {\scr D} = \scr D-\infty$ of degree $g-1$ in the definition of the Cauchy kernel, and then a Cauchy kernel with the above properties does not exist. For example the third property would become that $\div_q(C(p,q)) \geq - \wt {\scr D}-p$; now, $\wt {\scr D} +p$ is a generic divisor of degree $g$ (hence non-special) and then this forces $C(p,q)$ to be constant in $q$. Thus we cannot define a Weyl-Stieltjes differential and the corresponding Pad\'e\ problem.  See also Rem. \ref{remD2}.
 }

We now specialize the above general setup to the case that is of immediate interest in this work. 
\subsubsection{Real elliptic curves and asymptotics of orthogonal sections}
{  
 In general a smooth algebraic curve $\C$ is called ``real'' if there is an antiholomorphic map $\nu: \mathcal C\to \mathcal C$ which is involutive (i.e. $\nu^2 = Id$). In practice, if the curve is expressed in terms of a polynomial equation in two variables $W,X$, this means that we can choose the polynomial equation to have real coefficients and the antiholomorphic involution is typically $X\to \ov X, W\to \ov W$. 
In the case of an elliptic curve in Weierstrass form, like \eqref{Weierstrassparam} below,  it means that $g_2, g_3\in \R$.  Harnack's theorem states that the set of fixed points of $\nu$ consists of at most $g+1$ connected components (called ``real ovals'') and a curve with the maximum of real ovals is called a Harnack $M$--curve. The latter is the case we are considering here in genus $1$; in this case the modular parameter $\tau$ is shown to belong to $i\R_+$. For \eqref{Weierstrassparam} to be an $M$--curve, all roots must be real.} 

\begin{figure}
\begin{center}
\begin{tikzpicture}[scale=4]
\coordinate (tau) at (0,1.2);
\draw [fill=black!10!white] (0,0) to (tau)  to ($(tau)+ (1,0)$) to (1,0) to cycle;

\draw [blue,line width=1, postaction={decorate,decoration={{markings,mark=at position 0.7 with {\arrow[line width=1.5pt]{>}}}} }]($0.5*(tau)$) to node[pos=0,left] {$\frac \tau 2$} node[pos=0.6,below]{$\gamma$}($0.5*(tau)+(1,0)$);

\draw [blue,line width=1, postaction={decorate,decoration={{markings,mark=at position 0.7 with {\arrow[line width=1.5pt]{>}}}} }](0,0) to node[pos=0.6,below]{$\alpha$} (1,0);

\foreach \x in {0.0460022, 0.173841, 0.337538, 0.569893, 0.781358, 0.924899}
{\draw [fill, red] ($0.5*(tau)+ (\x,0)$) circle [radius=0.3pt];
}
\draw [fill, red] (0.4998,0) circle [radius=0.3pt];

\draw [fill, red!40!blue] (0.33333,0) circle [radius=0.31pt] node[below] {$\mathscr D$};
\node [below]at (0,0) {$0$};
\node [above]at (tau) {$\tau$};
\node [above]at ($(tau)+(1,0)$) {$\tau+1$};
\node [below]at (1,0) {$1$};
\end{tikzpicture}
\includegraphics{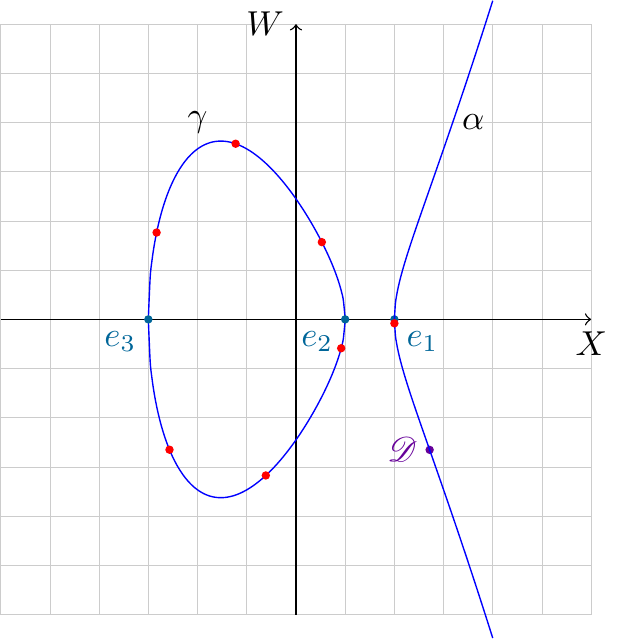}
\end{center}

\caption{
An example of real elliptic $M$-curve  (specifically $W^2 = 4(X-1)(X-2)(X+3)$). On the left pane we have the ``elliptic'' parametrization as the quotient of $\C$ by the lattice $\Lambda_\tau$. On the right the representation of the real section of $\mathcal E_\tau$  in the Weierstrass parametrization. The divisor $\scr D$ consists of a single point on the real oval of the $\alpha$ cycle (in this example $\scr D= 1/3$ in the elliptic parametrization), while the measure of orthogonality is defined on the cycle $\gamma$ and it is given by an arbitrary smooth positive weight ${\rm e}^{w(p)}$ on $\gamma$ times the holomorphic normalized differential $\d p = \frac {\d X}{2\omega_1 W}$. 
Also plotted are the zeros of the orthogonal section $\pi_6$ with respect to the ``flat'' measure with ${\rm e}^{w(p)} \equiv 1$. Note that the zero on $\alpha$ is already (for $n=6$) extremely close to $e_1$:  it is shown in Sec. \ref{epilogue} that  this zero, for even $n$ converges to $e_1$ exponentially fast.
}
\label{figone}
\end{figure}
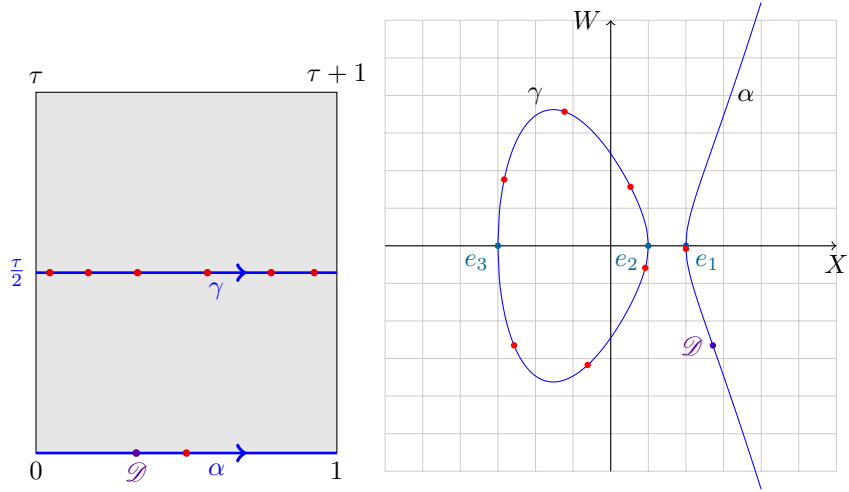

Thus, let now  $\mathcal E_\tau = \C/( \Z+\tau \Z)$ be a real  elliptic ($M$--)curve with $\tau \in i\R_+$. We can realize this curve equivalently as the locus of the algebraic equation in the $(X,W)$--plane
\be
\label{Weierstrassparam}
W^2 = 4 X^3 - g_2 X - g_3 = 4(X-e_1)(X-e_2)(X-e_3).
\ee
We refer to this representation as the ``Weierstrass parametrization''.
As mentioned above, the condition that this is an $M$--curve translates to the requirement that the roots are real and we order them as follows 
\be
e_3<e_2<e_1 \in \R.
\ee
This curve possesses the anti-involution given, in the two representations, by 
\be
(X,W)^\star = (\ov X, \ov W),\ \ \  \ p^\star = \ov p\in \C/(\Z+\tau \Z).
\ee
{ 
The two real ovals are simply the real components (i.e. with $X,W\in \R$) of the curve as shown in Fig. \ref{figone}. }
The relationship between the two representations is provided by the Weierstrass elliptic function $\wp$ as $X=\wp(p),\  W=\wp'(p)$.

The role of the point $\infty$ will be played by the point at infinity in the Weierstrass parametrization, or equivalently by $p=0$. 

There are two real components of the elliptic curve; the interval $[e_1, \infty)$ and $[e_3,e_2]$ on both sheets of the square--root $W = 2 \sqrt{(Z-e_1)(Z-e_2)(Z-e_3)}$, with the branch-cuts being chosen as the segments $(-\infty, e_3]\cup [e_2,e_1]$.

We will choose for $\gamma$ the real component $[e_3,e_2]$ (on both sheets); in the elliptic parametrization the contour $\gamma$ is represented by $\Im p=\frac {\Im \tau} 2$ and it is fixed by the anti-involution $^\star$ thanks to the fact that $\tau$ is purely imaginary.

The other fixed oval is $\alpha = [e_1,\infty)$ (on both sheets), which is also our ``alpha'' cycle. 
The $\beta$ cycle is the contour $[e_2,e_1]$ (again, on both sheets). They are represented in the elliptic representation as $\Im p\equiv 0\mod \Im \tau \Z$ and $\Re p\equiv 0\mod \Z$ respectively. See Fig. \ref{figone}.

Let $\d\mu$ be a positive measure on $\gamma$ of the form $\d\mu = {\rm e}^{w(p)}\d p$ with $w(p)$ a smooth real  function defined for $p\in \gamma$. We assume that the support of $\d\mu$ consists of the whole $\gamma$.  As customary we make the further assumption
\begin{assumption}
The function $w(p)$ is analytic in a strip containing $\gamma$ and real on $\gamma$. 
\end{assumption} 

To complete the data we need the divisor $\scr D$, which here consists of just one point. 
We choose it to be also invariant under the anti-involution $^\star$ and away from $\gamma$. Thus $\scr D$ consists of a single  point arbitrarily chosen on the $\alpha$--cycle (except $0$). In the elliptic representation it is chosen as $\scr D\in (0,1)$ modulo the lattice $\Lambda_\tau$. 

{  
The requirements that $\infty$ and $\scr D$ are fixed by the antiholomorphic involution is necessary to ensure that the space of meromorphic sections $\scr P_n = \mathcal O(n\infty +\scr D)$ has a real structure. This simply means that we can choose a basis that has the Schwarz property $f(p^\star) = \ov {f(p)}$.} 

With these preparations we are ready to formulate the results of this paper. 
\paragraph{Results.}
The first result is a characterization of the location of the zeros; namely we prove:\\[3pt]
\begin{theorem}[See Section. \ref{secexact}]  
The orthogonal sections $\pi_n$ exist and  have $n+1$ 
{ 
\em simple} zeros. These zeros lie all on  $\gamma$ for $n$ odd, while for $n$ even  one zero belongs to $\alpha$.
\end{theorem}

In Section \ref{secDZ} we tackle the much more intricate issue of the asymptotic behaviour. 
First of all the existence shown in Sec. \ref{secexact} implies that the solution of the RHP \ref{introYRHP} exists and is unique for every $n\in \mathbb N$. 

We then proceed to build a sequence of three transformations of the original problem into a final RHP \ref{RRHP}. This latter falls into the class of ``small norm'' Riemann--Hilbert problems.

The logic of these steps is completely parallel to the established literature. The main difference lies now in the analysis of the resulting integral equation (see Theorem \ref{smallR}). 

The obstacle is the issue of indices which we illustrate  here in an idealized situation.
Suppose we have a RHP for a matrix $Y(p)$ defined on the complement of the unit circle in the plane, and satisfying 
\be
Y_+(z) = Y_-(z)  J(z), \ \ \ |z|=1,\ \ \ Y(\infty)=\1.
\ee 
{
 Here the matrix--valued function $J: S^1\to GL_r(\C)$ is assumed to be  smoothly dependent on $z$. }
The solution we seek is assumed to be  such that both $Y(z), Y^{-1}(z)$ are bounded everywhere and analytic away from the circle. 
It is well known that the generic solvability of this problem requires that the index of $\det J$ is zero, namely, the phase of $\det J$ has null total increment around the circle.

One may naively assume that the same setup can be used on a Riemann--surface of higher genus by taking $\gamma$ to be the boundary of an embedded disk. However one quickly realizes \cite{Rodin}  that, 
{ 
for generic jump matrices $J$},   $\det Y$ necessarily must have $r\,g$ zeros and $r\, g$ poles (with $r$ being the size of the matrix) if $\det J$ has index zero. 
{  
To put it in a different, albeit possibly in a bit cryptic way, there are no holomorphic sections of a (non-trivial) vector bundle of degree zero on a Riemann surface of genus $g>0$; namely, the rows of $Y$ (which are the ``sections'' of the vector bundle) must have some poles. Then the index being zero implies that $\det Y$ has as many zeros as poles and the generic count is $r\,g$. 
 Rather than referring to \cite{Rodin}, which is a bibliographic reference that is a  very difficult to decipher for many readers (including the author of this manuscript), let us illustrate the obstacle by considering the simple example of $r=1$ (i.e. a line--bundle, or scalar Riemann Hilbert problem). 
Then it can be shown  that  the solution of the jump problem $Y(p_+)= Y(p_-) J(p)$ with $Y$ analytic and nonzero in $\mathcal C\setminus \gamma$ exists if and only if 
\be
\label{113}
\frac 1{2i\pi} \int_\gamma \ln J(p) \vec \omega(p) \equiv 0 \in \mathbb J(\mathcal C)
\ee
where $\vec \omega$ is the vector of normalized holomorphic differentials and $\mathbb J(\mathcal C)$ is the Jacobian of the curve. For example, in genus $1$ we have $\omega= \d p$ in the elliptic representation $p\in \C /(\Z+\tau \Z)$ and the condition reads 
\be
\label{114}
\frac 1{2i\pi} \int_\gamma \ln J(p) \d p =n + \tau m, \ \ \ n,m \in \Z.
\ee
We give a self-contained proof of this fact for genus $g=1$ in App. \ref{genex}; here it suffices to point out  that there are additional obstructions to the solvability of the RHP (even in the scalar case!) that did not exist in genus zero. 
As a remark for the reader oriented towards algebraic geometry we indicate  that the above is an explicit example of Picard's correspondence between line bundles of degree $0$ and points in the Jacobian;  the line bundle with transition function $J$ on the annulus around $\gamma$ is trivial if and only if the integral \eqref{114} (or \eqref{113} in higher genus) maps to the null-point in the Jacobian.\\[4pt]
}

We note that there is an obvious and inevitable connection with the general theory of vector bundles on Riemann  surfaces and this was extensively investigated in the recent \cite{BNR}.

The reader with some experience in the solution of the matrix problem may concretely appreciate the obstacle of solving a Riemann-Hilbert problem  by stepping through the initial phases of the solution; in genus zero the initial step consists in reformulating the problem as an integral equation of the form 
\be
Y(z) = \1  + \frac 1 {2i\pi} \oint_{|w|=1} Y_-(w) \big(J(w)-\1\big) \frac{\d w}{w-z}.
\ee
In attempting this step in higher genus, one needs to generalize the {\it Cauchy kernel}  $\frac {\d w}{w-z}$. Even for scalar kernels there is no unique generalization, and in fact the Cauchy kernels in higher genus depend on additional $g$ complex parameters, that manifest themselves in the choice of a divisor of degree $g$ \cite{Zvero, Fay}.
For matrix problems the issue is compounded and in fact the relevant notion is that of a {\it matrix} Cauchy kernel that depends on $r^2 g$ complex parameters ($r$ the size of $Y$) which go under the name of {\it Tyurin data}.  We refer to \cite{BNR} for a comprehensive description of these kernels. 

In our current setup we only see a faint shadow of all this rich theory: in genus $1$ much of the nuances can be disposed of because we have a global uniformization:
{  
 indeed, we can represent any elliptic curve as the quotient  $\C/(\Z+\tau \Z)$ and the coordinate on this plane (which we call $p$) provides a global coordinate; moreover the differential $\d p$ is  holomorphic and has no zeroes on the elliptic curve. This allows us  to think of differentials simply as elliptic functions by dividing them by $\d p$. This is problematic in higher genus since any  holomorphic differential has $2g-2$ zeros and hence dividing by any of them would introduce spurious poles. }
This simplification notwithstanding, we will encounter other novel features: notably,  in solving the small norm RHP \ref{RRHP}, and hence complete our steepest descent analysis, we will need a matrix (non-abelian) Cauchy kernel. This is described in Section \ref{Cauchykernel}.

\begin{figure}
\includegraphics[width=0.24\textwidth]{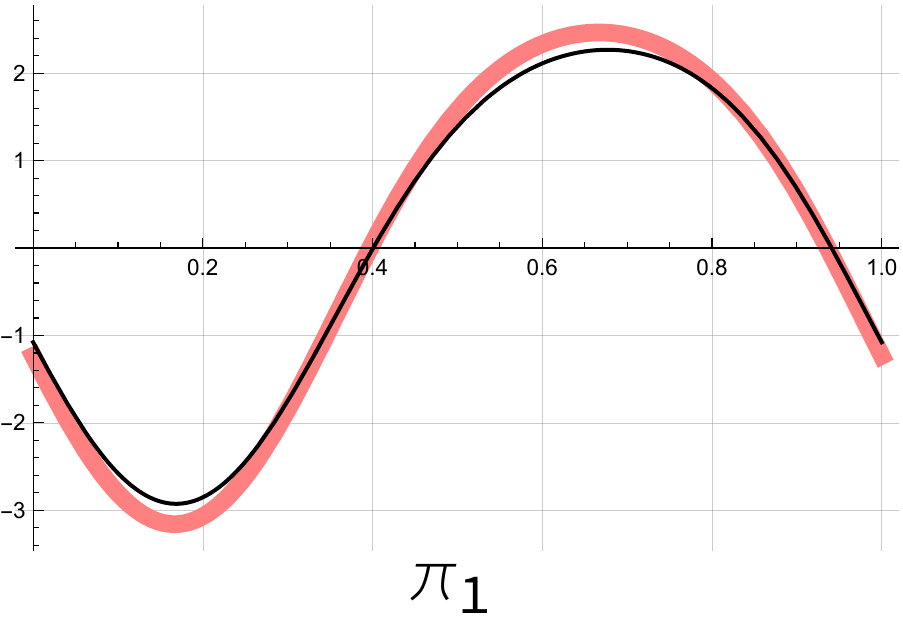}
\includegraphics[width=0.24\textwidth]{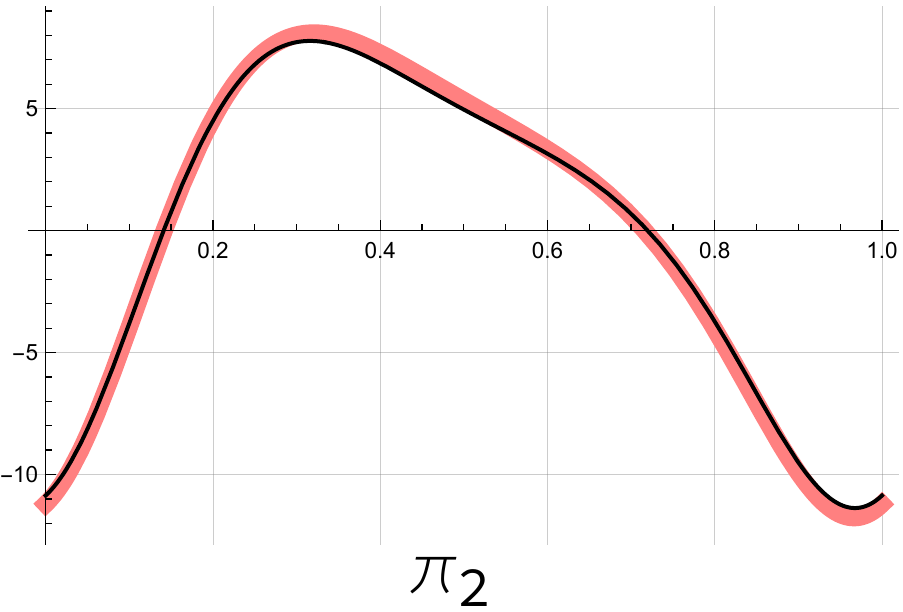}
\includegraphics[width=0.24\textwidth]{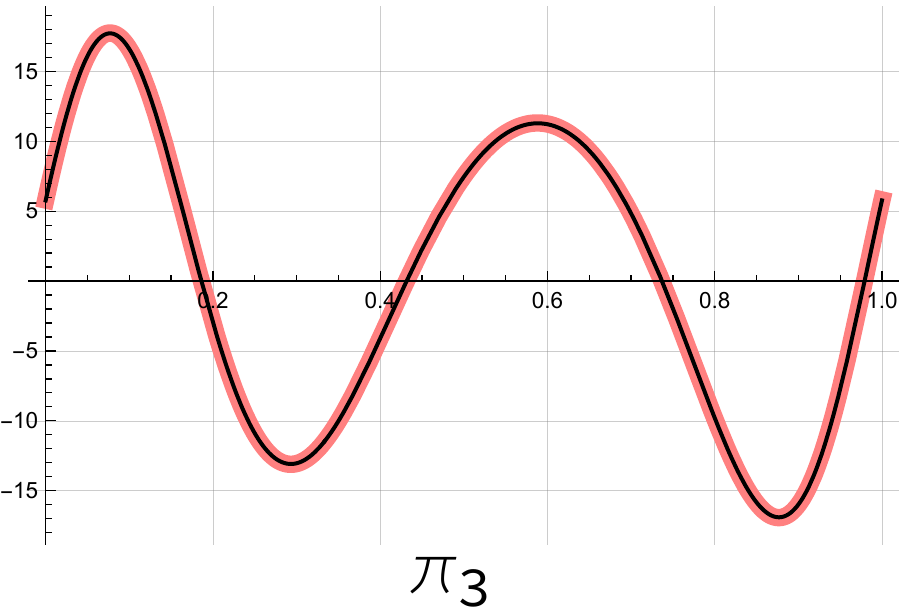}
\includegraphics[width=0.24\textwidth]{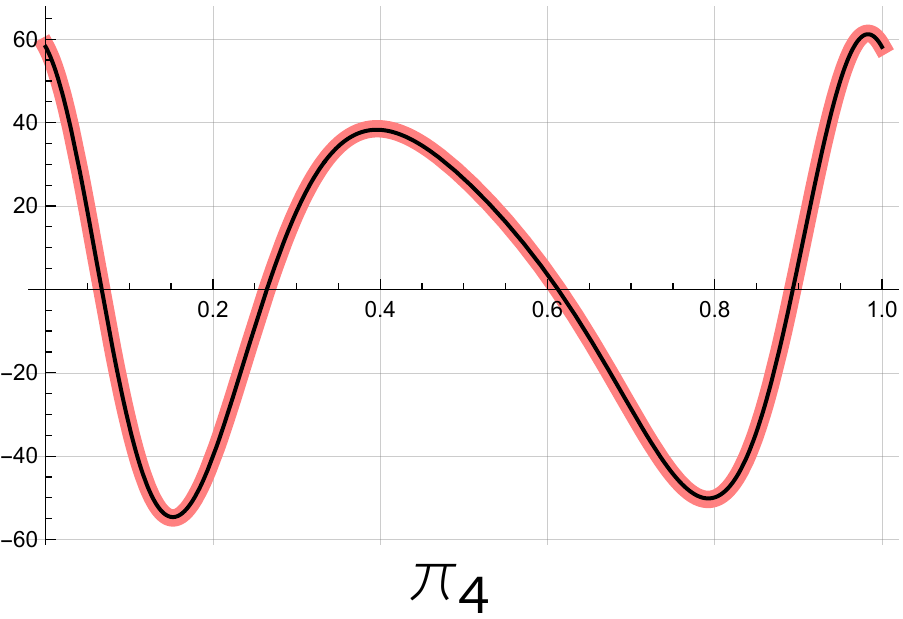}\\
\includegraphics[width=0.24\textwidth]{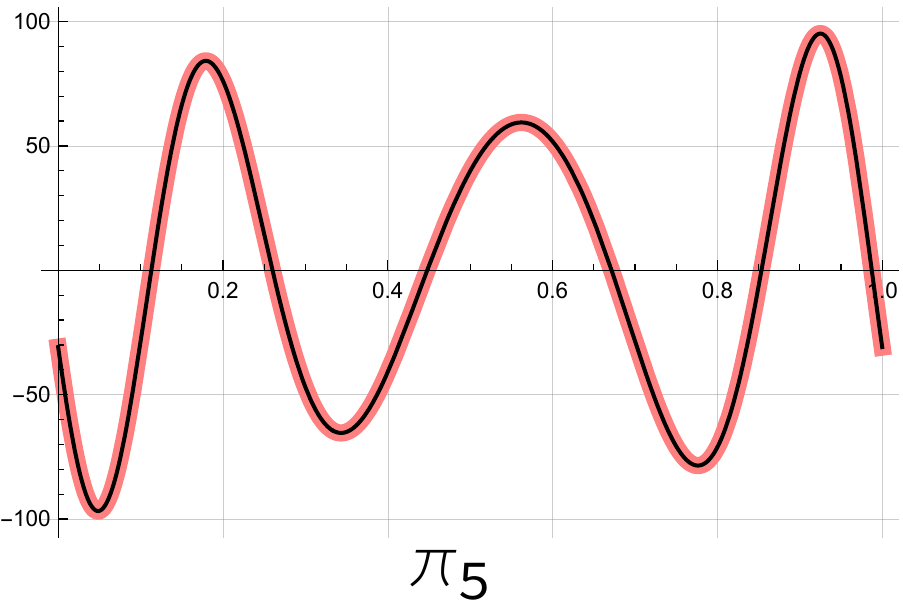}
\includegraphics[width=0.24\textwidth]{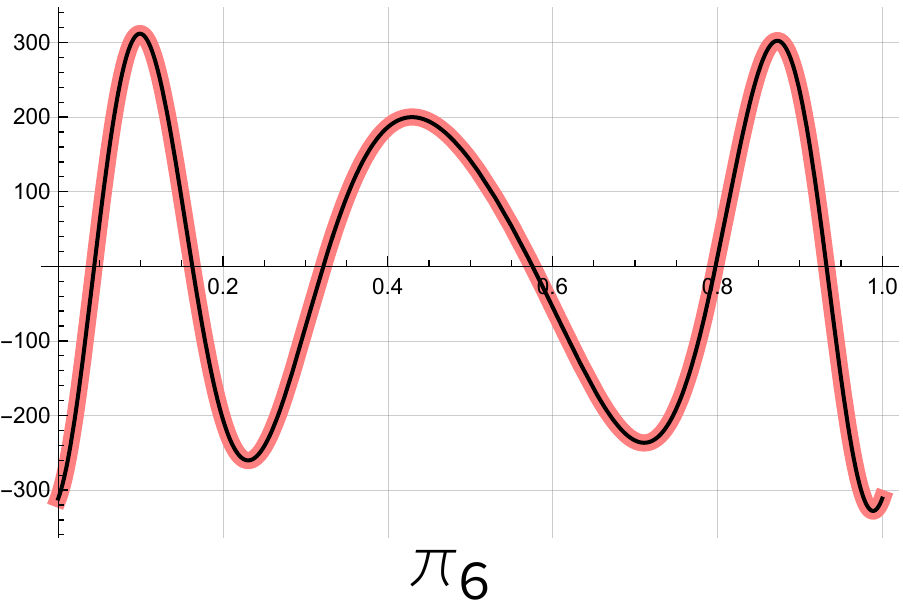}
\includegraphics[width=0.24\textwidth]{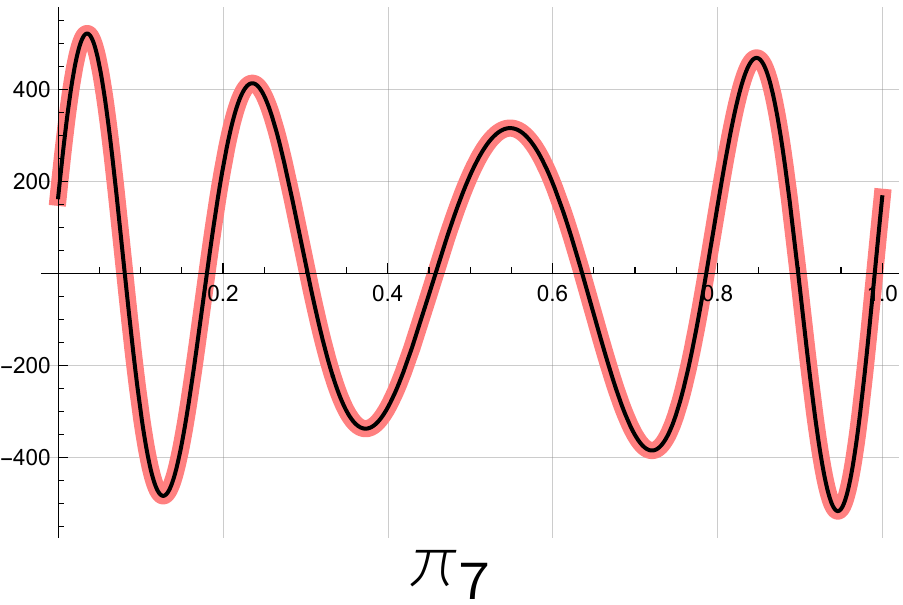}
\includegraphics[width=0.24\textwidth]{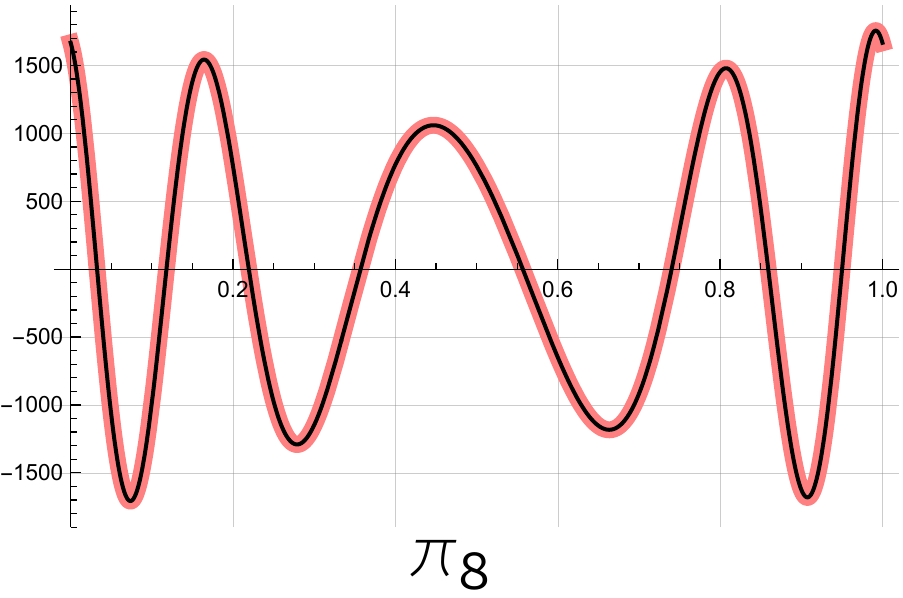}\\
\includegraphics[width=0.24\textwidth]{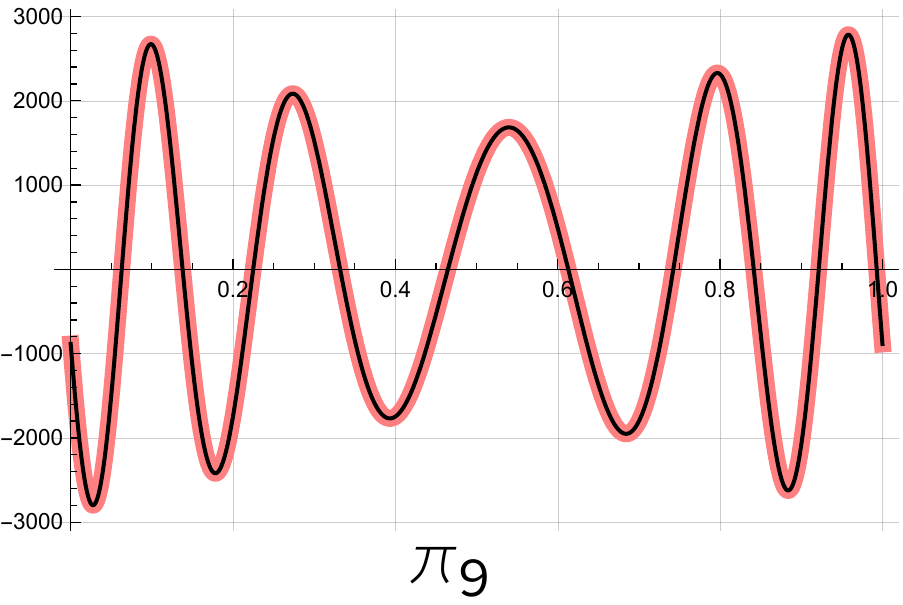}
\includegraphics[width=0.24\textwidth]{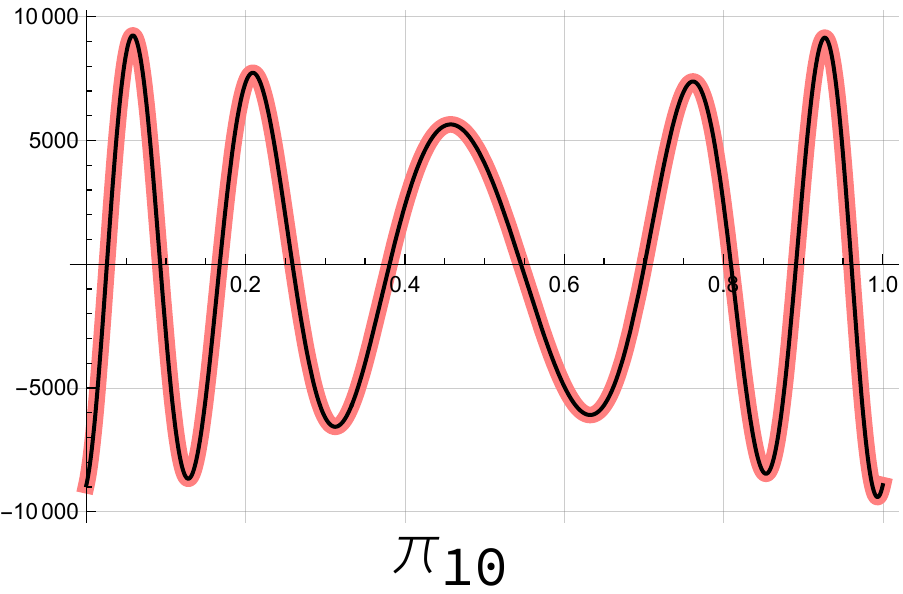}
\includegraphics[width=0.24\textwidth]{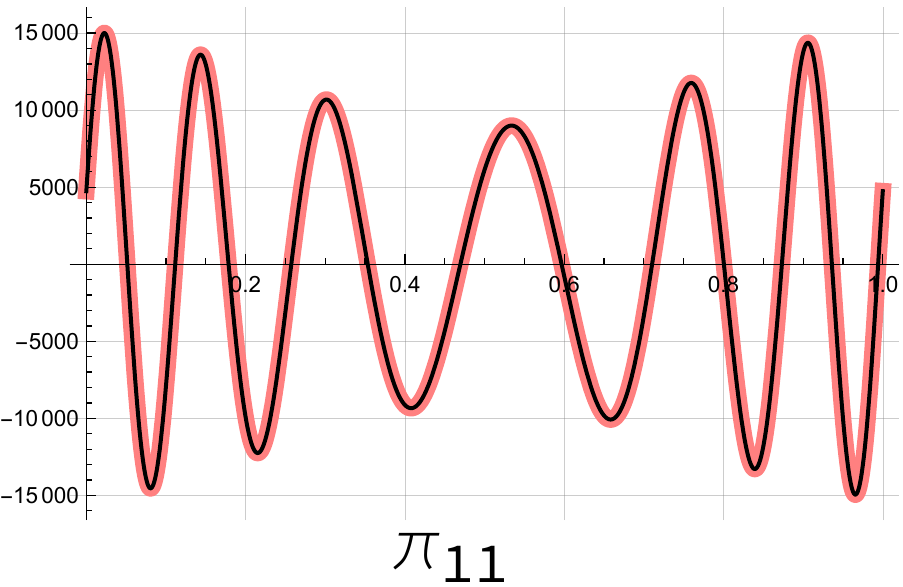}
\includegraphics[width=0.24\textwidth]{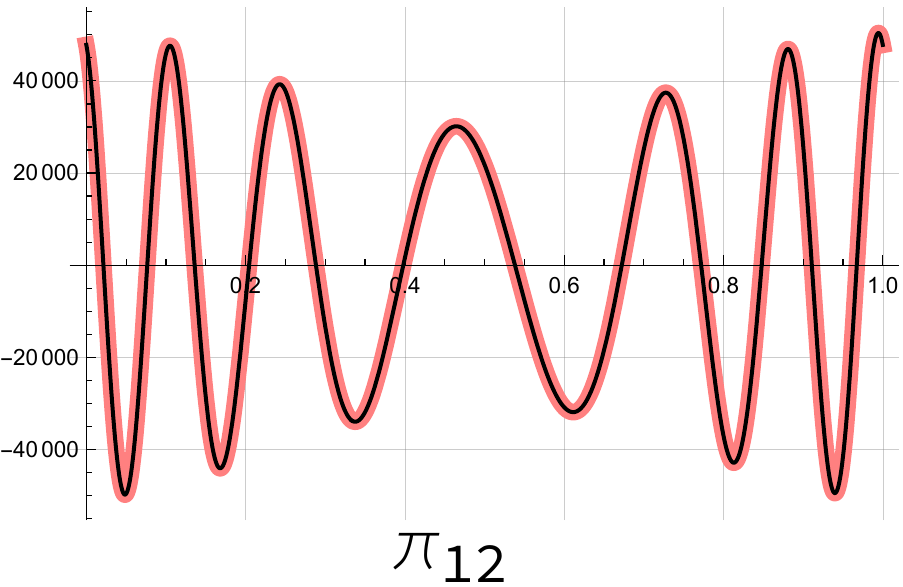}\\
\phantom{LA vispa teresa}

\caption{The first few monic orthogonal sections plotted as a function of $s\in [0,1]$ via $p=\frac \tau 2 +s$; here $\pi _n(p)\in \scr P_n$ are the ``monic'' sections behaving like $\pi_n(p) = p^{-n}(1 + \mathcal O(p))$. The elliptic curve is $W^2=4X^3 - 19X+15 = 4(X-1)(X-3/2)(X+5/2)$. Here $\tau \simeq 0.6563 i$. We have set $\scr D = 1/3\in \R$ and $\infty =0$. 
The contour $\gamma$ is  the segment $[\tau/2,\tau/2+1]$ in $\mathcal E_\tau$; in the $X$--plane this is the segment $X\in [e_3,e_2]$ (on both sheets).
The thick line is the plot of the orthogonal section obtained by computing explicitly the moments. The thin line is the approximation. Observe that the approximation is almost perfect starting from $n=2$, confirming the exponential rate of convergence discussed in the text.
 }
 \label{Figureplots}
\end{figure}

\begin{figure}
\centerline{\includegraphics[width=0.4\textwidth]{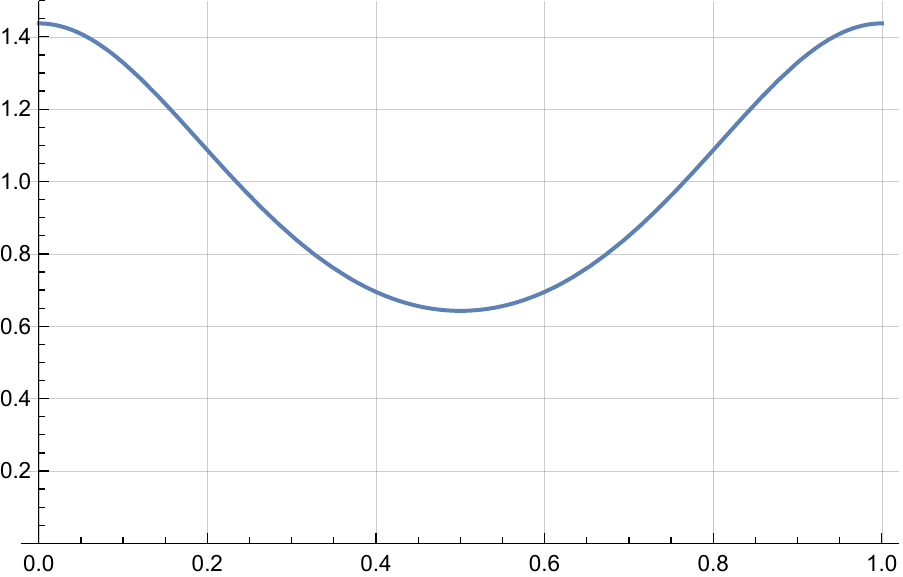}}

\caption{The asymptotic density $\d\mu_0/\d s$  of the counting measure of zeros of the orthogonal sections $\pi_n$ on $\gamma$, parametrized by $s\in [0,1]$. This corresponds to the elliptic curve $W^2 = 4(X-1)(X-2)(X+3)$.}\label{Figdensity}
\end{figure}

At the end of this analysis we have the following concrete characterization of the asymptotic behaviour (see Section \ref{epilogue}):
\begin{enumerate}
\item For every compact subset  of $\mathcal E_\tau \setminus \gamma$ we have
{  
\bea
\label{116}
\pi_n(p) =&   {\rm e}^{ - S_\infty}  \le( M_{11}(p)  +\frac{\mathcal O({\rm e}^{-n c_0})}{{\rm dist}( p, \scr D)} \ri)
{\rm e}^{(n-1) g(p) + S(p)}  
\eea}
The function 
$S(p)$ is the ``Szeg\"o'' function computed from the density ${\rm e}^{w(p)}$ in Definition  \ref{defSzego} (with $S_\infty=S(0)$). The Szeg\"o function is zero for the ``flat'' measure $\d\mu(p) = \d p$ (i.e. for $w(p)=0$).
The term 
  $M_{11}$ is given by  the following expression (see Theorem \ref{thmM})
\bea
M_{11}(p) = 
{\rm e}^{-i\pi p} \frac {\theta_1(\scr D ;2\tau)\theta_1(p-\scr D -\tau;2\tau)\theta_1'(0;2\tau) \theta_{\{2,3\}}(p;2\tau) }
{\theta_1(\scr D+\tau;2\tau)\theta_1(p-\scr D;2\tau)\theta_1(p;2\tau)\theta_{\{2,3\}}(0;2\tau)}.
\eea  
The choice between $\theta_2,\theta_3$ is made according to the parity of $n$; finally the $g$--function appearing in \eqref{116}  is given by the explicit expression (see Proposition  \ref{propg}):
\be
{\rm e}^{g(p)} ={\rm e}^\ell \le\{
\begin{array}{cc}
 \ds   {\rm e}^{i\pi \le( p -\frac \tau 2\ri) - \frac {i\pi}2 } \frac {\theta_1(p;2\tau)}{  \theta_1(p-\tau;2\tau)} & \Im \frac \tau 2< \Im p < \Im  \tau \\[14pt]
  \ds{\rm e}^{  -i\pi \le(p-\frac \tau 2\ri) +\frac {i\pi}2 } \frac{\theta_1(p-\tau;2\tau)} {\theta_1(p;2\tau)}  &0<\Im p< \frac 1 2 \Im \tau.
\end{array}
\ri.
\ \ \ \
{\rm e}^\ell= -i \frac { \theta_1'(0;2\tau)}{\theta_1(\tau;2\tau)} {\rm e}^{-i\pi \frac \tau 2}>0.
\ee
\item For $p\in \gamma$ we have the modulated oscillatory behaviour for $p = s+\frac \tau 2 +i0$:
{  
\be
\pi_n(p)=2 {\rm e}^{(n-1) \ell - S_\infty}  \le[  \Re \bigg(M_{11}(p_+) {\rm e}^{ S(p_+)}\le(  {\rm e}^{i\pi s - \frac {i\pi}2 } \frac {\theta_1\le(s+\frac \tau 2;2\tau\ri)}{  \theta_1\le(s-\frac \tau 2;2\tau\ri)}  \ri)^{n-1} \bigg)  + \mathcal O({\rm e}^{-nc_0})\ri]
\label{approxintrosupp}
\ee
}

The plots of the actual orthogonal sections (computed numerically) and the approximation \eqref{approxintrosupp} are shown together in Fig. \ref{Figureplots} for the flat measure $\d\mu(p)= \d p$, which clearly shows that the convergence is almost immediate.
\item Denoting by $z_{1}^{(n)},\dots, z_m^{(n)}$ with $m= 2\lfloor \frac {n+1} 2 \rfloor$ the zeroes of $\pi_n$ on $\gamma$ we have that for every continuous function $\phi$ defined on $\gamma\subset \mathcal E_\tau$ 
\be
\lim_{n\to\infty}\frac 1 n \sum_{j=1}^{ 2\lfloor \frac {n+1} 2 \rfloor} \phi(z_{j}^{(n)}) = 
\int_\gamma \phi(p)  \sqrt{e_1-\wp(p)} \frac{\d p}{2\pi}
\ee
See Fig. \ref{Figdensity} for a plot of this asymptotic density; also note that in the Weierstrass parametrization the asymptotic density is the usual  arcsine density on $[e_3,e_2]$ (divided by $2$ because $\gamma$ covers the interval twice).
\item The extra zero of $\pi_n$ for $n$ even  tends  at exponential rate to $p=\frac 1 2$ (i.e. $X=e_1$). 
\item The square of the norms of the monic orthogonal sections have the asymptotics (see Section \ref{norms})
\be
\|\pi_n\|^2 =
2\pi  {\rm e}^{2(n-1)\ell - 2S_\infty} {\rm e}^{-i\pi \tau} 
\frac { {\rm e}^{-2i\pi \scr D} \theta_1^2(\scr D ;2\tau)}{\theta_1^2(\scr D+\tau;2\tau)}
\frac {\theta_1'(0;2\tau) }{\theta_4(0;2\tau) }
 \le(\frac { \theta_3(0;2\tau)}{\theta_2(0;2\tau)  }\ri)^{(-1)^n} \big(1 + \mathcal O({\rm e}^{-nc_0})\big).
 \ee

\end{enumerate}

\paragraph{Notations.}
\begin{itemize}
\item[-] $\mathcal E_\tau$; the quotient of the plane $\C$ by the lattice $\Lambda_\tau = \Z+\tau \Z$, with $\Im \tau>0$. 
\item[-] $\theta_j(p;\tau)$ $j=1,2,3,4$ are  the four Jacobi theta functions. See \href{https://dlmf.nist.gov/20.2}{DLMF 20.2} or classical textbooks. We consider them as function of $p= \pi z$ in the DLMF normalization, so that, for example $\theta_1(p+1;\tau) = -\theta_1(p;\tau)$ rather than being $\pi$ quasi-periodic. 
\item[-] for an elliptic  function (or differential)  $f(p)$, $\div(f)$ denotes its divisor. To avoid confusion between the position of the point and the multiplicity, we surround the point by brackets. So, for example the divisor $2(1/3)+ (\tau/2)$ is the divisor consisting of the point $p=\frac 13$ with multiplicity $2$ and the point $p=\tau/2$ (modulo the lattice $\Lambda_\tau$).
\item[-] The big-oh notation $\mathcal O$ is used with two different meanings depending on the context. If $\scr D$  is any divisor, to say that $f\in \mathcal O(\scr D)$ means that $\div f\geq -\scr D$. In particular if $\scr D$ is positive then $f$ is allowed to have poles. If $\scr D$ is negative, $f$ must vanish at the points of $\scr D$. 

Viceversa, when we write $\mathcal O((p-a)^k)$ (for example) we understand this as the usual big-Oh notation in complex analysis, namely the germ of analytic  functions $f(p)$ such that $f(p)/(p-a)^k$ is locally bounded near $p=a$. The difference of usage should be clear by the context.
\item[-]
{  
 Similarly, $\mathcal K(\scr D)$ denotes the germs of differentials $\omega$ with pole/zeros as described above, namely such that $\div(\omega)\geq -\scr D$.}

\item[-] we use the normalized Weierstrass $\zeta, \wp$ functions (i.e. with a normalized lattice spanned by $1,\tau$) rather than the classical definition with a lattice spanned by $2\omega_1,2\omega_2$. 
Thus
\bea
\zeta(p) &= \frac 1 p + \sum_{n^2+m^2\neq 0} \le(
\frac 1{p + n + m\tau} - \frac 1{n + m\tau}  -\frac {p}{(n+m\tau)^2}
\ri),\\ 
\wp (p) &=-\zeta'(p) = \frac 1 {p^2} + \sum_{n^2+m^2\neq 0} \le(
\frac 1{(p + n + m\tau)^2}   -\frac {1}{(n+m\tau)^2}\ri).
\eea

\end{itemize}
\paragraph*{Acknowledgements.}
The work  was supported in part by the Natural Sciences and Engineering Research Council of Canada (NSERC) grant RGPIN-2016-06660.

\section{Existence of orthogonal sections and their zeros}
\label{secexact}
The existence was established already in \cite{BertoPade}, Section 2.3.1. To quickly review, we start by defining the initial basis of ``monic'' sections (relative to the choice of coordinate $p$) in $\scr P_n$, for example 
\be
\sigma_0=1, \ \ \sigma_1 = \zeta(p)-\zeta(p-\scr D), \ \ \  \sigma_{2\ell+2} = \wp^{\ell + 1}(p), \ \ \  \sigma_{2\ell+3} = -\frac 1 2\wp'(p) \wp^\ell (p), \ \ \ \ell \geq 0.
\ee
Note that only $\s_1$ has actually a pole at $\scr D$ while the other sections have only pole at the point $\infty$ (i.e. $p=0$).
\begin{remark}
In the Weierstrass parametrization \eqref{Weierstrassparam}, denoting the point $\scr D $ by the coordinate $(X_0,W_0)$, we can use (with a different normalization)
\be
\label{sigmas}
\sigma_0=1, \ \ \sigma_1 = \frac {W+ W_0}{X-X_0}, \ \ \sigma_{2\ell+2} = X^{\ell+1}, \ \ \sigma_{2\ell+3} = W X^\ell.
\ee
\end{remark}
We then  define the bi-moment matrix
\be
\mu_{ab} = \int_\gamma \sigma_a(p)\sigma_b(p){\rm e}^{w(p)}\d p.
\ee
Using the Andr\'eief identity for alternants \cite{Andreief}, one verifies that 
\be
D_n= \det \bigg[\mu_{ab}\bigg]_{a,b=0}^{n-1} = \frac 1 {n!} \int_{\gamma^n} \le(\det\bigg[\sigma_{a-1}(p_b)\bigg]_{a,b=1}^n\ri)^2 \prod_{a=1}^n{\rm e}^{ w(p_a)}\d p_a.
\ee
Since all functions $\s_\ell(p)$ are real--valued on $\gamma$ and ${\rm e}^{w(p)}\d p$ is a positive measure, it follows that $D_n>0$ and then the monic orthogonal sections are given by 
\be
\pi_n(p) = \frac 1{D_n} \det\le[
\begin{array}{cccccc}
\mu_{0,0} & \mu_{1,0} &\cdots & \mu_{n,0}\\
\mu_{0,1} & \mu_{1,1} &\cdots & \mu_{n,1}\\
\vdots &&&\vdots\\
\mu_{0,n-1} & \cdots && \mu_{n,n-1}\\
\sigma_0(p) & \sigma_1(p) &\cdots & \sigma_n(p)
\end{array}
\ri].
\ee

The non-vanishing of the determinant $D_n$ guarantees that the solution of the RHP \ref{introYRHP} (or equivalently the RHP \ref{YRHP} later) exists and is unique by Theorem 2.11 in \cite{BertoPade}.
 \par\vskip 4pt
 \begin{remark}
 \label{remD2}
 {  
We can see here that if $\scr D\to \infty = (0) $ (corresponding to $X_0\to \infty$ and hence also $W_0\to \infty$)  then  $\s_1$ in \eqref{sigmas} as written diverges. Of course we can normalize it by dividing by $W_0$ but then we obtain that $\s_1\to 1$ (uniformly over $\gamma$). This has the consequence that two rows of the matrix in the determinant $D_n$ are equal to each other and $D_n=0$. Consequently the RHP has no solution (as stated). It is not clear to the author if it is possible to consider a situation where $\infty$ belongs to $\scr D$ (even in higher genus). 
}
 \end{remark}

 We denote by $\scr P^\R$ the real span of the (real-analytic) sections $\{\s_\ell\}$. 
We now address the question of the positions of the zeros of the orthogonal sections $\pi_n$. 
\bet
\label{thmzeros}
{
The orthogonal sections $\{\pi_n\}_{n\in \mathbb N}$ have $n$ zeros on $\gamma$ and one on $\alpha$ for even $n$, while for odd $n$ it has $n +1$ zeros on $\gamma$ and none on $\alpha$.}
\eet
{\bf Proof.}
First of all note that $\pi_n$ has $n+1$ zeros since it has a pole of order $n$ at $p=0$ and of order $1$ at $\scr D$ and the total number of zeros of an elliptic function equals the total number of poles (counted with multiplicity). 

Since $\dim_\C\scr P_n=n+1$, we can construct a $\phi\in \scr P_n$ where we choose the positions of $n$ zeros while the remaining $g=1$ zero is uniquely determined by the Abel theorem, namely, that the sum of the positions of zeros must be congruent to the sum of the positions of poles modulo the lattice $\Lambda_\tau$.

Then we have that 
\be
\label{even}
\text{any real section in $\scr P^\R$ has an even number of zeros on $\gamma$, counting multiplicities.}
\ee
Indeed, let $\varphi\in \scr P^\R$; let $k$ be the number of zeros on $\gamma$, $\ell$ the number of zeros on $\alpha$ and $s$ the number of zeros elsewhere. 
By the Schwarz symmetry 
{ 
(namely, $f(\ov p) = \ov {f(p)}$)} it follows that $s = 2r$ is an even number.  The sum of the imaginary parts of these $2r$ zeros must be zero (modulo $\Im \tau$).  
Since the total of the sum of the { 
imaginary} parts { 
 of {\it all} the $\ell + s+k$ zeros} must vanish   (modulo $\Im \tau$) and $\Im \gamma = \Im \tau /2$, the statement follows by parity counting. 

Now suppose that $\pi_n$ has $2k$ zeros on  $\gamma$. If $2k\leq n-3$ then we can construct a real section $\varphi\in \scr P^\R_{2k+2}$ with $2k$ zeros on $\gamma$ and one zero on $\alpha$; it follows that the last zero (which we cannot choose arbitrarily) is nonetheless forced to lie in $\alpha$. 
Thus the product $\pi_n \varphi$ has constant sign on $\gamma$ and then, since the measure of orthogonality is positive, $\pi_n$  would not be orthogonal to $\varphi$, yielding a contradiction. We thus have determined that 
\be
\text{The orthogonal sections $\pi_n$ have at least $n-2$ zeros on $\gamma$. }
\label{atleast}
\ee
We now refine the result. Suppose first that $n$ is even; then $\pi_n$ has either $n$ or $n-2$ zeros on $\gamma$. Suppose the latter is the case. We can use a  real section of $\scr P_{n-1}$ (with a total of $n$ zeros) that has the same $n-2$ zeros on $\gamma$ and one (and hence both) remaining zeros on $\alpha$.  Thus the product $\pi_{n} \varphi$ has constant sign contradicting the orthogonality since $\int \varphi\pi_n {\rm e}^{w(p)}\d p$ cannot vanish. 

Suppose now that $n$ is odd; then $\pi_n$ has either $n+1$ or $n-1$ zeros on $\gamma$. Suppose that it has $n-1$ zeros: then, again, we can find $\varphi\in \scr P_{n-1}^\R$ with the same zeros on $\gamma$. The last zero is forced to lie on $\alpha$ by the Abel theorem. Proceeding as above we reach a contradiction. \QED

{  
\begin{corollary}
The zeros of $\pi_n$ are simple\footnote{We thank the anonymous referee for pointing out that the usual proof applies here.}.
\end{corollary}
{\bf Proof.} Let $z_1,\dots, z_k$ be the zeros of $\pi_n$ on $\gamma$ with odd multiplicity.
Note that $k$ must be even because the  total multiplicity of zeros on $\gamma$ is even.

Suppose first that $n$ is even, so that $\pi_n$ has a simple zero on $\alpha$ and $n$ zeros on $\gamma$ counted with multiplicity.   If $k<n$ (hence $k\leq n-2$ by parity) we can construct $\varphi\in \scr P^\R_{k}$ with simple zeros at $z_1,\dots, z_k\in \gamma$ and thus with a simple zero  $w_0\in \alpha$. The product $\varphi\pi_n$ does not change sign on $\gamma$ contradicting the orthogonality.

If $n$ is odd then all $n+1$ zeros of $\pi_n$ are on $\gamma$; again, the number of zeros with odd multiplicity is $k$ (even) and hence $k\leq n-1$. We can again construct  $\varphi\in \scr P^{\R}_k$ as above and reach the same conclusion.
\QED
}
\section{Steepest descent  analysis}
\label{secDZ}
In \cite {BertoPade} it was shown that the orthogonal sections $\phi_n$ are uniquely determined by a Riemann--Hilbert problem as in the RHP \ref{introYRHP}. 
Those details are a bit redundant in our present situation since the elliptic curve $\mathcal E_\tau$ admits a global uniformization and  moreover there is a differential without any zeros (i.e.  $\d p$) which we can use to trivialize the second column by simply dividing it by $\d p$. Thus we can reformulate the problem for a matrix consisting of bona-fide elliptic meromorphic functions with a jump discontinuity across $\gamma$.  We obtain the following RHP, which is thus the starting point of our investigation.
\begin{shaded}
\begin{problem}
\label{YRHP}
Let $Y= Y_n(p)$ be the $2\times 2$ matrix, meromorphic on $\mathcal E_\tau \setminus \gamma$ and with  poles at $p=0, \scr D$, such that 
\begin{enumerate}
\item
Near $p=0\equiv \Lambda_\tau$ we have the behaviour
\be
Y(p) = \le(\1 + \mathcal O(p) \ri) \le[
\begin{array}{cc}
p^{-n}  & 0 \\
 0 & p^{n-2}
\end{array}
\ri],\ \  \ p\to 0 \mod \Lambda_\tau.
\ee
\item Near $p=\scr D \mod \Lambda_\tau$ we have that 
\be
Y(p) = \le[\begin{array}{cc}
\mathcal O( (p-\scr D)^{-1})   &\mathcal O(p-\scr D )\\
\mathcal O( (p-\scr D)^{-1})   &\mathcal O(p-\scr D )
\end{array}\ri],
\ee
namely, the first column has at most a simple pole and the second column has at least a simple zero at $p=\scr D$ (we use here the complex-analysis convention since $p$ is a global coordinate) . 
\item The boundary values at $p\in \gamma$ are bounded and satisfy:
\be
Y(p_+) = Y(p_-) \le[
\begin{array}{cc}
1 & {\rm e}^{w(p)}\\
0 & 1
\end{array}
\ri].
\ee
\end{enumerate}
\end{problem}
\end{shaded}
We note that $\det Y(p)$ is a meromorphic elliptic function without a jump across $\gamma$; it has a double pole at $p=0\mod \Lambda_\tau$ and hence it must have two zeros in each fundamental domain. These points are the {\it Tyurin points} of a vector bundle that depend on $n$ and the orthogonality measure. 
\subsection{The three  transformations}
The general approach of the steepest descent method requires the construction of a suitable $g$--function, as the nomenclature goes.

In our case, since the measure of orthogonality is independent of $n$, the $g$--function is a universal object independent of the particular measure as long as  its support is the whole $\gamma$ (which is our standing assumption).

Consider the expression 
\be
\h(p) = \sqrt{\wp (p)- \wp(1/2)} = \sqrt{X-e_1}.
\ee
The expression in the square-root has a single double zero at $p=\frac 1 2\mod \Lambda_\tau$. Thus the square-root can be defined in a neighbourhood of the zero in $\mathcal E_\tau$ as a single valued function. However the expression is certainly not a perfect square because, if that were the case, the result would be  an elliptic function with a single simple pole at $p=0$ and a single simple zero at $p=\frac 1 2$. Rather, classical formul\ae\ (see \href{https://dlmf.nist.gov/23.6.E5}{DLMF 23.6.5}) show that this extends to a multi-valued  elliptic function with the property 
\be
\h(p+1) = \h(p), \ \ \ \h(p+\tau)= -\h(p).
\ee
We take avail of this multivaluedness on $\mathcal E_\tau$ by defining:
\be
\label{defH}
H(p):= \le\{
\begin{array}{cc}
-\h(p) & 0< \Im(p) < \frac \tau 2\\
\h(p) & \frac \tau 2< \Im(p)<\tau
\end{array}
\ri..
\ee
Observe that $H(p)$ now has a discontinuity across $\gamma$  and satisfies $H_+ + H_-=0$,  but it is single--valued  on the remainder of the elliptic curve, namely $H(\tau+ s) = H(s)$, $s\in \R$. It has a simple pole at $p=0\mod \Lambda_\tau$ and 
\be
\label{resH}
\res{p=0} H(p) \d p = -1.
\ee
In terms of the Weierstrass parametrization we have 
\be
\label{3.8}
H(p) \d p =  \sqrt{X-e_1}\frac {\d X} {2 \sqrt{(X-e_1)(X-e_2)(X-e_3)}}  = 
\frac {\d X} {2 \sqrt{(X-e_2)(X-e_3)}},
\ee
where the determination is such that the radical behaves like $X$ at $X=\infty$.   We define now the $g$--function.
\begin{definition}
\label{defGfun}
The $g$--function is defined by the following integral 
\be
g(p) =\lim_{t\to 0} \le( \ln t + \int _t^p H(q)\d q\ri).
\ee
\end{definition}
\bp
\label{propg}
The $g$--function has the properties; 
\begin{enumerate}
\item $g(p) =- \ln p + \mathcal O(p)$, as $p\to 0$ and has two branch-cuts along the $\gamma$ and $\alpha$ cycle;
\item $g(p+1)= g(p) -  i\pi$ for $\frac 1 2 \Im \tau < \Im p <\Im \tau$ and
$g(p+1)= g(p) +  i\pi$ for $0<\Im p<\frac 1 2 \Im \tau $. In particular along the $\alpha$ cycle it has $g(p_+)= g(p_-)+2i\pi$; 
\item $\Im g(p_+)=-\Im g(p_-)$ is a monotonically decreasing function  for $\Im p= \frac 12 \Im \tau$ as $\Re p$ increases. Moreover $\Im g(p_\pm) = \mp \pi \int_{\frac \tau 2 } ^p \d\mu_0(q)$, with $\d\mu_0 = \sqrt{e_1 - \wp(p)} \frac{\d p}{\pi}$ a smooth probability measure on $\gamma$. 
\item The $g$ function is given explicitly by the following expression 
\be
{\rm e}^{g(p)} =K \le\{
\begin{array}{cc}
 \ds   {\rm e}^{i\pi \le( p -\frac \tau 2\ri) - \frac {i\pi}2 } \frac {\theta_1(p;2\tau)}{  \theta_1(p-\tau;2\tau)} & \Im \frac \tau 2< \Im p < \Im  \tau \\[14pt]
  \ds{\rm e}^{  -i\pi \le(p-\frac \tau 2\ri) +\frac {i\pi}2 } \frac{\theta_1(p-\tau;2\tau)} {\theta_1(p;2\tau)}  &0<\Im p< \frac 1 2 \Im \tau.
\end{array}
\ri.
\ \ \ \
K= -i \frac { \theta_1'(0;2\tau)}{\theta_1(\tau;2\tau)} {\rm e}^{-i\pi \frac \tau 2}
\label{explicitG}
\ee
and it satisfies
\be
{\rm e}^{g_+(p)+ g_-(p)} =K^2  = -\frac { \theta_1'(0;2\tau)^2}{\theta_1(\tau;2\tau)^2} {\rm e}^{-i\pi \frac \tau 2} =: {\rm e}^{2\ell}, \ \ \ p\in \gamma.
\label{g+-}
\ee
where $\ell\in \R$.
\item Alternatively, in the Weierstrass  parametrization \eqref{Weierstrassparam} we have the simple expressions 
\bea
\label{gsimple}
{\rm e}^{g(p)- \ell} =\frac{2}{e_3-e_2}\le[ \le(X-\frac {e_2+e_3}2\ri) +  \sqrt{(X-e_2)(X-e_3)}\ri],\ \ \ X = \wp(p),
\eea
where the branch-cut of the radical is on $[e_3,e_2]$ with the radical behaving as $X$ for $X\to\infty$. This expression has modulus $1$ for $X_\pm \in [e_3,e_2]$ (on both sides) and of modulus strictly greater than $1$ on $\C \setminus [e_3,e_2]$.
\end{enumerate}
\ep
{\bf Proof.}
{\bf (1.)} The first property follows from the fact that the differential $\d g = H(p)\d p$ has residue $-1$ at $p=0$. The absence of the constant is due to the definition of the regularization. 
\\
{\bf (2.)} 
By construction/definition 
\be
\d g(p_+) + \d g(p_-)=0, \ \ p\in \gamma, \  \ \res{p=0} \d g=-1.
\ee
Integrating the jump condition along $\gamma$ (which is homologous to the $\alpha$--cycle) we get
\be
\int_{\gamma + i\epsilon}  \d g + \int_{\gamma-i\epsilon} \d g=0.
\ee
Now we use Cauchy's residue theorem to ``wrap around'' the contour $\gamma - i\epsilon$ by increasing $\epsilon=0 \to \Im \tau$; since $\d g$ has a single pole at $p=0$ on the torus, we pick up the minus the residue which gives a  contribution of $2i\pi$. Thus we get 
\bea
0=\int_{\gamma + i\epsilon}  \d g + \int_{\gamma-i\epsilon} \d g= \int_{\gamma + i\epsilon}  \d g + \int_{\gamma+i\epsilon-\tau} \d g + 2i\pi   = 2\int_{\gamma + i\epsilon} \d g +  2i\pi.
\eea
The last equation implies $g(p+1)-g(p)= \int_p^{p+1} \d g =- i\pi$, for $\Im\tau>\Im p>\frac 1 2\Im \tau$. The other case is treated similarly.
\\
{\bf (3.)} By construction $\d g(p_\pm) = \mp \sqrt{X-e_1} \d p$ with $X = \wp (p)$. Since $X\in [e_3,e_2]$ this expression is imaginary and clearly $\d g(p_+)=-\d  g (p_-)$ for $p\in \gamma$. Now, by an application of the Cauchy residue theorem, we have 
\be
-2i\pi= 2i\pi \res{p=0} \d g = \int_{\gamma} \d g(p_+)-\d g(p_-)= 2\int_\gamma \d g(p_+).
\ee
Hence the total mass of $\d\mu_0 = \frac {\sqrt{e_1-X}}{\pi} \d p$ is one.\\
{\bf (4.)} 
The function $\h(p) = \sqrt{\wp (p)-e_1}$ extends to an elliptic function on $\mathcal E_{2\tau}$ with simple poles at $p=0,\tau$ modulo $\Lambda_{2\tau}$ and residues $\pm 1$, respectively. Thus one can write it as $\h(p) = \zeta(p;2\tau) - \zeta(p+\tau;2\tau) + c$ for a suitable constant. Given that $\zeta$ is minus  the logarithmic derivative of $\theta_1$ up to an additive constant, one deduces that the antiderivative of $\h(p)$ must be of the form 
\be
\int \h(p)\d p = \ln \frac {\theta_1(p+\tau;2\tau)}{\theta_1(p;2\tau)} + c p + d
\ee
with $c, d$ constants. The constant $c$ is determined by the fact that $g(p+1) = g(p)- i\pi$  and hence ${\rm e}^{g(p)}$ must be anti-periodic  for $p\to p+1$. The constant $d$ of integration is determined by the regularization that requires that ${\rm e}^{g(p)} = \frac 1 p(1+\mathcal O(p))$ near $p=0$. Then the reader may verify that the expression \eqref{explicitG} satisfies all these requirements. Also, one should verify that ${\rm e}^{g(p)}$  is continuous on the torus away from $\gamma$. The rest is an exercise using the known properties of the Jacobi theta functions. 
 The constant $K$ is simply the normalization so that ${\rm e}^{g(p)} = p(1 + \mathcal O(p))$ near $p=0$.  Finally the fact that $\ell \in \R$ is due to the fact that the expression 
\be
{\rm e}^{2\ell} = -{\rm e}^{-i\pi \tau}\frac{ \theta_1'(0;2\tau)^2} {\theta_1^2(\tau;2\tau)}>0 
\ee 
is positive thanks to $\theta_1(\tau;2\tau)$ being purely imaginary.\\
{\bf (5.)} We know that $g(p)-\ell$ is purely imaginary on the cycle $\gamma$ and hence the modulus of ${\rm e}^{g(p)-\ell}$ is one on the $\gamma$ cycle. Furthermore we know that 
\be
\d g(p) = \sqrt{\wp(p) - e_1} \d p =\frac 1 \pi  \sqrt{X -e_1} \frac {\d X}{\sqrt{4(X-e_1)(X-e_2)(X-e_3)}}
= \frac {\d X}{2\pi\sqrt{(X-e_2)(X-e_3)}}.
\ee
Note that the integration on $\gamma$ translates, in Weierstrass parametrization, to the segment $[e_3,e_2]$ covered {\it twice}. From elementary integration one obtains the claim. 
Finally one observes that \eqref{gsimple} is precisely the conformal map of the plane $\C$ slit across $[e_2,e_3]$ to the outside of the unit disk; from this we deduce that the modulus is $>1$ outside of the cycle $\gamma$. 
\QED

We will also need the ``Szeg\"o" function; for that purpose we need the following (scalar) Cauchy kernel. This is a differential $C(p,q) \d p$  where $C(p,q)$ here we think of  as a meromorphic elliptic function of the two variables $p,q$:
\be
C(p,q ) = \frac {\theta_1'(0) \theta_1(q) \theta_1(p-1/2) \theta_1 (p-q+1/2)}
{\theta_1(1/2) \theta_1(p) \theta_1(p-q) \theta_1(q-1/2)}.
\ee
It has the defining properties
\begin{enumerate}
\item $\div_q(C(p,q)) \geq -p + (0) - (1/2)$;
\item $\div_p(C(p,q)) \geq -q - (0) + (1/2)$;
\item $\res{p=q} C(p,q)\d p=1$. 
\end{enumerate}
In the Weierstrass parametrization, denoting the two points by the coordinates $p\equiv(X,W), \ q\equiv (X',W')$, it has the simple expression
\be
C(p,q)\d p = \le(\frac {W+W'}{X-X'} - \frac { W'}{e_1-X'} \ri) \frac {\d X}{2W}.
\ee
We can thus proceed to the definition of the Szeg\"o function.
\begin{definition}
\label{defSzego}
 The Szeg\"o function for the measure density ${\rm e}^{w(p)}$ is defined by:
\be
S(p) = \int_\gamma C(p, q)  w(q) \frac {H(q_+) \d q}{H(p)2i\pi} = \int_\gamma C(p, q) w(q) \frac {\d g(q_+)}{2i\pi\frac {\d g(p)}{\d p}} 
\ee
where $H(p) = \frac {\d g(p)}{\d p}$ is as in \eqref{defH}.
\end{definition}

The Szeg\"o function has the following easily established properties 
\begin{enumerate}
\item $S(p_+) + S(p_-) =-  w(p), \ \ \ p\in \gamma$;
\item $S(p+1)=S(p)$,  $ S(p+\tau) = S(p)$;
\item $S(p)$ is otherwise a bounded analytic function on $\mathcal E_\tau \setminus \gamma$;
\item It satisfies the Schwarz symmetry $\ov S(\ov p) =  S(p)$, $p\in \mathcal E_\tau \setminus \gamma$, and for $p\in \gamma$ we thus have $\ov{S(p_+)} = S(p_-)$; 
\item The Schwarz symmetry, together with the fact that $w(p)$ is a real--valued function on $\gamma$, imply that 
\be
\label{SSchwartz}
S(p_\pm ) = \pm i\,\nu(p) - \frac 1 2 w(p),\ \ \ p\in\gamma.
\ee
where $\nu(p)$ is a real--analytic function defined in a small strip around $\gamma$ and $\nu(p+1)=\nu(p)$. 
\end{enumerate}

Finally we introduce the {\it effective potential}
\be
\label{defphi}
 \varphi(p;n)= \varphi(p) := (n-1)  \big(g(p)-\ell\big)+ S(p) +\frac 12  w(p).
\ee
which is analytic within the strip of analyticity of $w(p)$ minus $\gamma$. 
The following Lemma illustrates its main properties.
\begin{lemma}[Properties of the effective potential]
\label{lemmaphi}
The effective potential satisfies:
\begin{enumerate}
\item $\varphi(p_+)+ \varphi(p_-) = 0$ for $p\in \gamma$.
\item Suppose that the domain of analyticity of $w(p)$ contains the $\epsilon_0$ strip around $\gamma$. Then, for every $0<\epsilon<\epsilon_0$  there is $N_0$ such that if $n>N_0$ the real part of $\varphi$ on  $\gamma_\pm := \gamma \pm  i\epsilon$ is strictly positive. Furthermore 
\be
\inf_{p\in \gamma+i\epsilon} \Re \varphi(p;n) \mathop{\to}_{n\to\infty} +\infty\ .
\ee
\item We have $\varphi(p_\pm) = \mp i\pi (n-1)  \int_{\frac \tau 2}^p \sqrt{e_1-\wp(q)} \d q  \pm i\nu (p)$, where $\nu$ is the phase function defined in \eqref{SSchwartz}.
\end{enumerate}
\end{lemma}
{\bf Proof.} The property {\bf (1.)} is simply a rearrangement of \eqref{g+-} and the properties of the Szeg\"o\ function.
\\
{\bf (2.)} We have shown  in point $3$ of Proposition \ref{propg} that $\Im g(p_+)$  is a monotonically strictly decreasing  function. By the Cauchy--Riemann relations applied to $g(p)$ from the $+$ side, we deduce that the real part is increasing as we move upwards. Similar considerations apply to the negative boundary value, which has the property that $\Im g(p_-) = - \Im g(p_+)$  is {\it increasing}, and hence the real part increases also moving downwards. Thus $\Re (g(p)-\ell)$ is strictly positive in a strip  around $\gamma$ (in fact, we have proved that it is strictly positive everywhere out of $\gamma$).  In general we cannot say much about the real part of $S(p)+ \frac {w(p)} 2$ within its analyticity domain. However, since $g(p)-\ell$ is multiplied by $(n-1) $ in the definition of $\varphi$ \eqref{defphi}, for every $\epsilon$ we clearly have that the real part of $(n-1) (g(p)-\ell)$ eventually dominates over that of $S(p) + \frac {w(p)} 2 $ and also tends to infinity.
\\
{\bf (3.)} Since $\varphi_++\varphi_-=0$ and also $\varphi(\ov p) = \ov{\varphi(p)}$ we conclude that $\varphi(p_+) = -\varphi(p_-) \in i\R$. 
Then, from \eqref{SSchwartz} we obtain.
\be
\varphi(p_\pm) = \mp i\pi (n-1)  \int_{\frac \tau 2}^p \sqrt{e_1-\wp(q)} \d q  \pm  i\nu (p)\ .
\ee
Note that $\varphi(p+1)=\varphi(p)\mp i\pi$ in the two halves of $\mathcal E_\tau$ bounded by the $\gamma$ and $\alpha$ cycles.  
\QED

With the definitions of the $g$  and Szeg\"o\ functions we can now start our chain of transformations. 
\subsubsection{First transformation}
Let 
\be
\label{Zdef}
Z(p) ={\rm e}^{(- (n-1)  \ell + S_\infty)\s_3} Y(p) {\rm e}^{-((n-1)  (g(p)-\ell)  + S(p))\s_3}\ .
\ee
We then obtain a new Riemann--Hilbert problem for the matrix $Z$, which is a simple exercise starting from the properties of $Y(p)$ in RHP \ref{YRHP} the definition \eqref{defphi} of $\varphi$ and its properties in Lemma \ref{lemmaphi}.
\begin{shaded}
\begin{problem}[The RHP for $Z$]
\label{ZRHP}
The matrix--valued function $Z(p)$ { 
is meromorphic  in $\E_\tau \setminus \gamma$ and it}  satisfies:
\begin{enumerate}
\item
Near $p=0\equiv \Lambda_\tau$ we have the behaviour
\be
Z(p) =  \le[
\begin{array}{cc}
p^{-1}  + \mathcal O(1) & \mathcal O(1) \\
 \mathcal O(1) & p^{-1} + \mathcal O(1)
\end{array}
\ri],\ \  \ p\to 0 \mod \Lambda_\tau\ .
\ee
\item Near $p=\scr D \mod \Lambda_\tau$ we have that 
\be
Z(p) = \le[\begin{array}{cc}
\mathcal O( (p-\scr D)^{-1})   &\mathcal O(p-\scr D )\\
\mathcal O( (p-\scr D)^{-1})   &\mathcal O(p-\scr D )
\end{array}\ri]\ .
\ee
\item The boundary values at $p\in \gamma$ are bounded and satisfy:
\bea
Z(p_+) &= Z(p_-) \le[
\begin{array}{cc}
{\rm e}^{-(n-1)  \Delta g- \Delta S }& {\rm e}^{(n-1) (g_+ + g_- -2\ell)+ S_+ +  S_-+w }\\
0 & {\rm e}^{(n-1)  \Delta g + \Delta S }
\end{array}
\ri]
\\
=&Z(p_-) \le[
\begin{array}{cc}
1& 0 \\
{\rm e}^{-2\varphi_- } & 1
\end{array}
\ri]
 \le[
\begin{array}{cc}
0 & 1\\
-1 & 0
\end{array}
\ri]
 \le[
\begin{array}{cc}
1& 0 \\
{\rm e}^{-2\varphi_+ } & 1
\end{array}
\ri]\ .
\eea
where, for a function $F$ with a discontinuity across $\gamma$ we have denoted $\Delta F(p) = F(p_+)-F(p_-)$. 
\item We have the periodicity $Z(p+1)=(-1)^{(n-1) } Z(p)$.
{  
Moreover, for $p\not \in \gamma$, we have $Z(p+\tau)= Z(p)$. }
\end{enumerate} 

\end{problem}
\end{shaded}

\subsubsection{Opening lenses: second transformation.}
Let us fix $\epsilon>0$ (small) such that $w(p)$ is analytic in the whole open  $\epsilon$ strip around $\gamma$. Consider then the matrix
\be
\label{Tdef}
T(p)= \le\{
\begin{array}{cc}
Z(p) & \le|\Im (p-\frac \tau 2)\ri| >\epsilon \\
Z(p)  \le[
\begin{array}{cc}
1& 0 \\
-{\rm e}^{-2\varphi_+ } & 1
\end{array}
\ri]
 &  0<\Im (p-\frac \tau 2) <\epsilon
\\
Z(p)  \le[
\begin{array}{cc}
1& 0 \\
{\rm e}^{-2\varphi_- } & 1
\end{array}
\ri]
 &  -\epsilon <\Im (p-\frac \tau 2) <0
\end{array}
\ri.\ .
\ee
Once more a direct verification shows that this new matrix satisfies the following problem.
\begin{shaded} 
\begin{problem}[$T$--problem]
\label{TRHP}
The matrix--valued function $T(p)$ is meromorphic in $\E_\tau\setminus \gamma$ and it satisfies:
\begin{enumerate}
\item
Near $p=0\equiv \Lambda_\tau$ we have the behaviour
\be
T(p) = \le[
\begin{array}{cc}
p^{-1}  + \mathcal O(1) & \mathcal O(1) \\
 \mathcal O(1) & p^{-1} + \mathcal O(1)
\end{array}
\ri],\ \  \ p\to 0 \mod \Lambda_\tau\ .
\ee
\item Near $p=\scr D \mod \Lambda_\tau$ we have that 
\be
T(p) = \le[\begin{array}{cc}
\mathcal O( (p-\scr D)^{-1})   &\mathcal O(p-\scr D )\\
\mathcal O( (p-\scr D)^{-1})   &\mathcal O(p-\scr D )
\end{array}\ri]\ .
\ee
\item on the two contours $\gamma_\pm := \gamma \pm  i\epsilon$ we have 
\be
T(p_+)  = T(p_-) \le[
\begin{array}{cc}
1& 0 \\
{\rm e}^{-2\varphi(p)  } & 1
\end{array}
\ri], \ \ \ p\in \gamma_\pm \ .
\label{lensjump}
\ee
\item On the contour $\gamma$ we have:
\bea
T(p_+) &= T(p_-) \le[
\begin{array}{cc}
0 & 1\\
-1 & 0 
\end{array}
\ri]\ .
\eea
\item We have the periodicity $T(p+1)=(-1)^{(n-1) } T(p)$, $T(p+\tau)= T(p)$ for $p$ outside of $\gamma$.
\end{enumerate} 
\end{problem}
\end{shaded}

It follows now from the property 2. in Lemma \ref{lemmaphi} that  the two matrices featuring in the jumps relations \eqref{lensjump} are (uniformly) exponentially close to the identity matrix as $n\to\infty$. More precisely there is a positive constant 
\be
\label{defc0}
c_0 = \frac 1 2  \min_{p\in \gamma_\pm} \Re (g(p) -\ell)>0
\ee
 such that $\le|{\rm e}^{-2\varphi(p)} \ri| <  {\rm e}^{-n c_0}$ for $n$ sufficiently large.
 
To complete the asymptotic analysis we have to construct a ``model solution'' (see next section) which gives an exact solution of a problem  similar to the RHP \ref{ZRHP} but where \eqref{lensjump} is absent. We then show  that the remainder terms are small; this logic follows the usual thread of ideas as in \cite{DKMVZ} but there are differences due to the fact that the model problem has a determinant that vanishes at certain points ({\it Tyurin points}).

\subsubsection{Model problem}
The model parametrix must satisfy  the following problem
\begin{shaded}
\begin{problem}[$M$-problem]
\label{MRHP}
The matrix--valued function $M(p)$ is meromorphic in $\E_\tau \setminus \gamma$ and it satisfies:
\begin{enumerate}
\item
Near $p=0\equiv \Lambda_\tau$ we have the behaviour
\be
M(p) =  \le[
\begin{array}{cc}
p^{-1}  + \mathcal O(1) & \mathcal O(1) \\
 \mathcal O(1) & p^{-1} + \mathcal O(1)
\end{array}
\ri],\ \  \ p\to 0 \mod \Lambda_\tau\ .
\label{M0}
\ee
\item Near $p=\scr D \mod \Lambda_\tau$ it satisfies the growth behaviour 
\be
\label{Mdiv}
M(p) = \le[\begin{array}{cc}
\mathcal O( (p-\scr D)^{-1})   &\mathcal O(p-\scr D )\\
\mathcal O( (p-\scr D)^{-1})   &\mathcal O(p-\scr D )
\end{array}\ri]\ .
\ee
\item On the contour $\gamma$ the boundary values are related by:
\bea
M(p_+) &= M(p_-) \le[
\begin{array}{cc}
0 & 1\\
-1 & 0 
\end{array}
\ri]\ .
\label{Mjump}
\eea
\item The matrix $M$ has the periodicity $M(p+1)=(-1)^{(n-1) } M(p)$,  $M(p+\tau) = M(p)$ for $p\not\in \gamma$.
\end{enumerate} 
\end{problem}
\end{shaded}
\begin{figure}
\begin{center}
\begin{tikzpicture}[scale=4]
\pgfmathsetmacro{\R}{4} ;
\pgfmathsetmacro{\RR}{int(\R+\R)} ;
\pgfmathsetmacro{\Rd}{int(\R-1)} ;
\coordinate (tau) at (0,1.2);
\draw [fill=black!10!white] (0,0) to (tau)  to ($(tau)+ (1,0)$) to (1,0) to cycle;

\fill [red!10!white] 
($0.6*(tau)$) to ($0.6*(tau)+(1,0)$) to  ($0.5*(tau)+(1,0)$) to ($0.5*(tau)$);
\fill [pink] 
($0.4*(tau)$) to ($0.4*(tau)+(1,0)$) to  ($0.5*(tau)+(1,0)$) to ($0.5*(tau)$);

\draw [green!50!black,line width=1, postaction={decorate,decoration={{markings,mark=at position 0.7 with {\arrow[line width=1.5pt]{>}}}} }]($0.5*(tau)$) to node[pos=0,left] {$\frac \tau 2$} node[pos=0.6,below]{$\gamma$}($0.5*(tau)+(1,0)$);

\draw [blue!50!black,line width=1, postaction={decorate,decoration={{markings,mark=at position 0.7 with {\arrow[line width=1.5pt]{>}}}} }]($0.6*(tau)$) to node[pos=0.4,above]{\tiny $\gamma_+$}($0.6*(tau)+(1,0)$);

\draw [blue!50!black,line width=1, postaction={decorate,decoration={{markings,mark=at position 0.7 with {\arrow[line width=1.5pt]{>}}}} }]($0.4*(tau)$) to node[pos=0.4,below]{\tiny $\gamma_-$}($0.4*(tau)+(1,0)$);

\node at ($0.5*(tau) + (0.8,0.05)$) {$\mathcal L_+$};
\node at ($0.5*(tau) + (0.2,-0.06)$) {$\mathcal L_-$};

\node [below]at (0,0) {$0$};
\node [above]at (tau) {$\tau$};
\node [above]at ($(tau)+(1,0)$) {$\tau+1$};
\node [below]at (1,0) {$1$};
\draw [fill, red!40!blue] (0.33333,0) circle [radius=0.31pt] node[below] {$\mathscr D$};

\end{tikzpicture}
\end{center}
\caption{The strip-lens for the RHP \ref{RRHP}. 
}
\label{figlenses}
\end{figure}
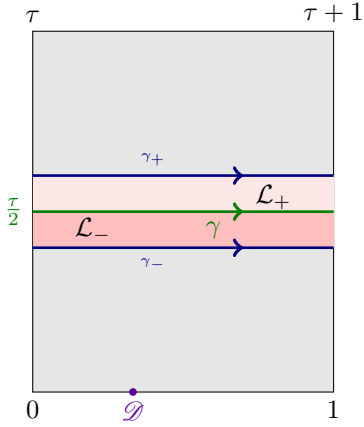
The following theorem and subsequent corollary show that the solution exists and is unique.

\bet
\label{thmM}
The solution to the RHP \ref{MRHP} exists and is unique; it is given by
\be
M(p) =C\le\{\begin{array}{lc}
\ds
 \le[\begin{array}{cc}
f(p) & f(p+\tau)\\
h(p) & h(p+\tau)
\end{array}\ri] & 0 \leq \Im p<\frac 1 2 \Im \tau
\\[10pt]
\ds
 \le[\begin{array}{cc}
-f(p + \tau) & f(p)\\
-h(p+\tau) & h(p)
\end{array}\ri] &\frac 1 2 \Im \tau< \Im p<\Im \tau,
\end{array}
\ri.
\ee
where $f, h$ are the two following functions on the elliptic curve $\mathcal E_{2\tau}$
\be
f(p):=
\le\{
\begin{array}{cc}
\ds 
{\rm e}^{-i\pi p} \frac {\theta_1(p-\scr D -\tau;2\tau) \theta_2(p;2\tau) }
{\theta_1(p-\scr D;2\tau)\theta_1(p;2\tau)}
 & \text{$ n  $ even}\\[14pt]
\ds
{\rm e}^{-i\pi p} \frac {\theta_1(p-\scr D -\tau;2\tau) \theta_3(p;2\tau) }
{\theta_1(p-\scr D;2\tau)\theta_1(p;2\tau)}
 & \text{$n $ odd }
\end{array}\ri.
,\ \ \ \ 
\label{defF}
\ee
\be
h(p):= 
\le\{
\begin{array}{cc}
\ds 
{\rm e}^{-i\pi p+ i\pi \tau + 2i\pi \scr D- \frac {i\pi}2 } \frac {\theta_1(p-\scr D -\tau;2\tau) \theta_3(p;2\tau) }
{\theta_1(p-\scr D;2\tau)\theta_4(p;2\tau)} 
& \text{ $n $ even}\\[14pt]
\ds 
{\rm e}^{-i\pi p+ i\pi \tau + 2i\pi \scr D-\frac {i\pi}2}  \frac {\theta_1(p-\scr D - \tau;2\tau) \theta_2(p;2\tau) }
{\theta_1(p-\scr D;2\tau)\theta_4(p;2\tau)}
 & \text{ $n $ odd }
\end{array}\ri. .
\label{defG}
\ee
The functions $f,h$ also have the following properties on the doubled elliptic curve $\mathcal E_{2\tau}$, which characterize them up to a multiplicative scalar:
\bea
&\div_p(f) \geq -(0)  -(\scr D) + (\scr D+\tau)
\qquad f(p+1)=(-)^{(n-1) }f(p);  \ \ \ \ \ f(p+2\tau) = - f(p) .
\label{propf}
\\
& \div_p(h)\geq  -(\tau) -(\scr D) + (\scr D+\tau) 
\qquad h(p+1)= (-)^{(n-1) } h(p),\ \ \ h(p+2\tau) =- h(p).
\label{defg}
\eea
The normalization matrix $C$ is 
\be
\label{defC1}
C =\le(\frac { \theta_1(\scr D ;2\tau)}{\theta_1(\scr D+\tau;2\tau)}\ri)^{\s_3}
\frac {\theta_1'(0;2\tau)}{\theta_2(0;2\tau) },\ \ \ \s_3:= \le[\begin{array}{cc}
1 & 0\\
0&-1
\end{array}\ri] 
\ee
for $n$ even and 
\be
\label{defC2}
C =\le(\frac { \theta_1(\scr D;2\tau)}{\theta_1(\scr D+\tau;2\tau)}\ri)^{\s_3}
\frac {\theta_1'(0;2\tau)}{\theta_3(0;2\tau) }\ee
for $n $ odd.
\eet

\noindent{\bf Proof.} The existence is simply a verification of the indicated formula; the uniqueness follows from the argument below.

Let $f(p)= M_{11}(p)$;  it follows from the  Problem \ref{MRHP} and in particular \eqref{Mjump} that $f(p)$ admits analytic continuation across $\gamma$ and that the resulting function (denoted by the same symbol) is a function with the properties \eqref{propf} on the double covering elliptic curve $\mathcal E_{2\tau}$. Similar considerations applied to $h(p)= M_{21}(p)$ show that it extends to a meromorphic function with the properties \eqref{defg}.
{  
The minus sign  $f(p+2\tau) = -f(p)$ in \eqref{propf} (and similarly in \eqref{defg}) is due to the fact that continuation across $\gamma$ (see \eqref{Mjump}) {\it twice} (i.e. going around the $b$--cycle of $\E_\tau$ twice) brings back $f$ to $- f$ (simply because the square of the constant jump matrix in \eqref{Mjump} is minus the identity). }

These properties determine uniquely the functions $f,h$ from a simple exercise in the use of Jacobi theta functions.

Finally, the normalization behaviour \eqref{M0} is then guaranteed by the matrix $C$. Note that the zeros of $\theta_{1,2,3,4}(p;2\tau)$ inside the fundamental domain are at $0,\frac 1 2, \tau+\frac 1 2, \tau$, respectively. Therefore the matrix $C$ is well defined (i.e. the denominators are non-zero).
\QED
\begin{corollary}
\label{cornonabtheta}
Consider the RHP \ref{MRHP} where the condition \eqref{M0} is replaced by the condition of being locally holomorphic at $p=0$. Then the only solution is the trivial one. 
\end{corollary}
{\bf Proof.}
Suppose that $M(p)$ is a non-trivial solution of the stated problem (i.e. analytic at $p=0\in \mathcal E_\tau$). We can follow the same proof as for Theorem \ref{thmM} and define similarly the function $f=M_{11}$ and $h=M_{21}$. They will now have the divisor properties 
\be
\div_p(f)= \div_p(h) =-(\scr D)+ (\scr D+\tau)\ ,
\ee
(since they have a single pole they can have only one zero) 
and have the periodicities  indicated in \eqref{propf}, \eqref{defg}.

If we integrate $\frac{p \d \ln f}{2i\pi} $ along the boundary of the fundamental domain we obtain the sum of the positions of the zeros minus the position of the poles. However, the quasi--periodicity of $f$ allow us to deduce 
\be
\tau = \oint_{\pa \mathcal E_{2\tau}} \frac {p f'(p)}{f(p)}\frac{ \d p}{2i\pi} = - \int_{\pa \mathcal E_{2\tau}} \ln f(p) \frac{\d p}{2i\pi}  = \le\{
\begin{array}{cc}
\frac 1 2  + \tau\mod \Z &\text {if $(n-1) $ is odd}\\
\frac 1 2 \mod \Z &\text {if $(n-1) $ is even\ ,}
\end{array}
\ri.
\ee
which gives a contradiction. \QED

A second consequence of the uniqueness is in the following corollary.
\begin{corollary}
\label{SchartzM}
The matrix $M(p)$ enjoys the symmetries:
\be
\ov M(\ov p) =  i^{\s_3} M(p) i^{-\s_3}.
\ee
In particular $\ov M_{11}(\ov p) = M_{11}(p)$  and, for $p\in \gamma$ we have  $M_{12}(p_+) = M_{11}(p_-) = \ov{ M_{11}(p_+)}.$
\end{corollary}
{\bf Proof.} As usual one observes that the matrix $\wt M(p):= \ov {M(\ov p)}$   satisfies the same RHP \ref{MRHP}  except for the condition \eqref{Mjump} which is replaced by 
\be
\wt M(p_+) = \wt M(p_-) \le[
\begin{array}{cc}
0 & -1\\ 1 &0
\end{array}
\ri]
\label{Mtildejump}
\ee
because the conjugation interchanges the sides of the boundary values. 
Observing that the jump matrix in \eqref{Mtildejump} is the conjugation of the one in \eqref{Mjump} by $i^{\s_3}$, we arrive at the conclusion  that the matrix $i^{-\s_3} \wt M(p) i^{\s_3}$ satisfies the same RHP \ref{MRHP}. By  the established uniqueness of the solution in Theorem \ref{thmM} we conclude that $\wt M$ must coincide with $M$.\QED 

\subsubsection{Tyurin data}
We now begin to see a first main difference from the usual setup; indeed it is immediate to realize that the determinant of $M$ has exactly two zeros. 
This follows from the fact that $\det M(p)$ is an elliptic function with a double pole at $p=0$ and hence it must have two zeros. By the Schwarz symmetry $\det M(p)$ is a real--analytic elliptic function and hence the zeros form an invariant set under the anti-involution.

We now find the zeros of $\det M$ by direct inspection.
\begin{proposition}[Tyurin data]
\label{proptyurin}
The determinant of $M(p)$ vanishes at the {\it Tyurin divisor} $\scr T = (1/4) + (3/4)$. The spanning vectors of the left null-spaces of $M(p)$ at these two points are linearly independent.
\end{proposition}
{\bf Proof.}
Consider the case of $(n-1) $ odd. The matrix $M$ can be factored as 
\be
M = C \le[
\begin{array}{cc}
\frac{ \theta_2(p;2\tau)}{ \theta_1(p;2\tau)} & 
\frac{ \theta_2(p+\tau ;2\tau)}{ \theta_1(p+\tau;2\tau)} \\
\frac{ \theta_3(p;2\tau)}{ \theta_4(p;2\tau)} & 
\frac{ \theta_3(p+\tau ;2\tau)}{ \theta_4(p+\tau;2\tau)} 
\end{array}
\ri] {\rm diag} \le({\rm e}^{-i\pi p } \frac {\theta_1(p-\scr D - \tau;2\tau)}{\theta_1(p-\scr D;2\tau)},{\rm e}^{-i\pi p-i\pi \tau } \frac {\theta_1(p-\scr D ;2\tau)}{\theta_1(p-\scr D + \tau;2\tau)}
 \ri)
\ee 
and it is easily seen that the determinant of the leading and trailing matrices cannot vanish: actually the determinant of the last term is ${\rm e}^{-2i\pi \scr D-i\pi \tau}$ and it is constant with respect to $p$. The middle matrix can be rewritten as 
\be
\le[
\begin{array}{cc}
\frac{ \theta_2(p;2\tau)}{ \theta_1(p;2\tau)} & 
\frac{i \theta_3(p;2\tau)}{ \theta_4(p;2\tau)} \\
\frac{ \theta_3(p;2\tau)}{ \theta_4(p;2\tau)} & 
\frac{ i\theta_2(p ;2\tau)}{ \theta_1(p;2\tau)} 
\end{array}
\ri]\ .
\ee
The determinant can be simplified using the addition formula \href{https://dlmf.nist.gov/20.7.E7}{DLMF 20.7.7} to give 
\be
\det M \propto \frac{ \theta_2^2(p;2\tau)\theta_4^2(p;2\tau)- \theta_3^2(p;2\tau)\theta_1^2(p;2\tau)}{\theta_1^2(p;2\tau)\theta_4^2(p;2\tau)}
\propto \frac{ \theta_2(2p;2\tau)}{\theta_1^2(p;2\tau)\theta_4^2(p;2\tau)},
\ee
where the proportionality is up to constants independent of $p$. The numerator vanishes at $2p = \frac 1 2 + \Z + 2\tau \Z$ and hence at the points $p=\frac {1}4 , \frac 34$ in $\mathcal E_\tau$. The case of even $(n-1) $ is the same. 

To show finally that the two left null spaces at the points of the Tyurin divisor $\scr T$ are linearly independent we could find them explicitly or, more elegantly, reason as follows. 
Suppose by contrapositive  that ${\bf h}$ is  a nontrivial row--vector spanning simultaneously the two left null spaces of $M(p)$ at both points. 
Then  consider the following vector--valued function:
\be
{\bf f}(p):= \omega (p) {\bf h} M(p), \  \ \ \ \ \omega(p):= \frac {\theta_1^2(p)}{\theta_1(p-\frac 1 4)\theta_1(p-3/4)}\ .
\ee
Here $\omega$ is an elliptic function with two simple poles at $\scr T$ and one double zero at $p=0$.   The result is that ${\bf f}(p)$ is analytic at $\scr T$ because ${\bf h} M(\scr T)=0$ and vanishes as $p=0$.   But this would contradict the already established uniqueness of the solution of the RHP \ref{MRHP}, because we could add ${\bf f}(p)$ to a row of $M$ and obtain another solution of the same RHP. \QED

\subsubsection{Third transformation}
We now denote by $R(p)$ the matrix
\be
R(p) = T(p) M(p)^{-1}\label{Rdef}.
\ee
 Direct analysis shows that it satisfies the following problem;
\begin{shaded}
\begin{problem}[Remainder term problem]
\label{RRHP}
The matrix $R(p)$ is meromorphic on $\mathcal E_\tau\setminus (\gamma_+\cup \gamma_-)$, with bounded boundary values on $\gamma_\pm$  satisfying:
\be
R(p_+) = R(p_-) M(p) \le[\begin{array}{cc}
1 & 0\\
{\rm e}^{-2\varphi(p)} & 1
\end{array}\ri] M^{-1}(p), \ \ \ \ p\in \gamma_\pm,
\ee
where $\varphi$ is  the effective potential \eqref{defphi}.
The only singularities of $R$ are 
{
simple poles}  at the Tyurin divisor $\scr T$ (the zeros of the determinant of $M(p)$) and they are such that the product $R(p) M(p)$ is locally analytic. 
Finally, the matrix $R$ is normalized by the condition $R(0)=\1$
\end{problem}
\end{shaded}

\subsection{Non-abelian Cauchy kernel for vector bundles and small-norm theorem}
\label{Cauchykernel}
The sequence of transformations leading to the Riemann--Hilbert problems \ref{ZRHP}, \ref{TRHP} and ultimately to \ref{RRHP} follows closely the general strategy in the established literature on the ``Deift--Zhou'' method \cite{DKMVZ}. The last step is however significantly different because it relies upon the solvability (for large $n$) of the RHP for the remainder term. This problem is affected by the topology of the Riemann surface as discussed in the introduction. To proceed we need to define a non-abelian (i.e. matrix) version of the ordinary Cauchy kernel. The data that define this kernel are the Tyurin data mentioned in the introduction and in the previous section.
 
\paragraph{Definition of the non-abelian Cauchy kernel.}
The non-abelian Cauchy kernel $\Cauchy_0(p,q)\d p$ \cite{BNR} is a $2\times 2$ matrix--valued  differential with respect to the variable $p$ and meromorphic matrix--valued function with respect to the variable $q$ satisfying the following properties
\begin{enumerate}
\item It has a simple pole for $p=q$ and $p=0$ and no other poles with respect to $p$;
\item The residue matrix for $p=q$ is $\1$ (and hence at $p=0$ is $-\1$) 
\item It has a simple pole for $q=p$ and at the Tyurin divisor $\scr T= (1/4) + (3/4)$ and all entries vanish for $q=0$. 
\item The expression $M^{-1}(p) \Cauchy_0(p,q) M(q)$ is locally analytic with respect to $q$ and $p$ at $\scr T$. 
\end{enumerate}
These conditions determine the kernel uniquely as we now show. An explicit expression is given in \cite{BNR} in arbitrary genus,  but it will not feature prominently here.

\begin{remark}
The Cauchy kernel depends only on the Tyurin divisor $\scr T$ and on the null-spaces of $M$ at the Tyurin points. See \eqref{Cexp} and following.
\end{remark}

\paragraph{Uniqueness.}
Suppose we have another ``Cauchy'' kernel $\mathbf F(p,q)\d p$ satisfying the same conditions. 
{ 
Let us choose $C = M(p)^{-1}$ for  $p\not\in \scr T$ and}  consider the two matrices  $S_1(q):=C \mathbf F(p,q)M(q)$, $S_2(q):= C\Cauchy_0(p,q) M(q)$ viewed as  functions of $q$. They satisfy the same  RHP \ref{MRHP} save for the fact that the pole behaviour \eqref{M0} is not at $q=0$  but at $q=p$.  
{ 
Indeed, both  Cauchy kernels have the local behaviour $\frac {\1}{p-q}  + \mathcal O(1)$ as $p\to q$.}
Then the difference $S_{12}(q):= S_1(q)-S_2(q)$  has no pole at $q=p$ and is therefore a solution of the RHP \eqref{MRHP} which is also analytic at $q=0$. Then $S_{12}$ is necessarily zero as was shown in Corollary \ref{cornonabtheta}.

\paragraph{Explicit expression (not necessary).}
Let ${\bf v}_1, {\bf v}_2$ be two spanning row-vectors of the left null-spaces  of $M(p)$ at the two Tyurin points $\scr T = (1/4)+(3/4) = t_1+t_2$.

Then the Cauchy kernel is given by 
\bea
\label{Cexp}
\Cauchy_0(p,q)\d p  = \omega_{0,q}(p)
\1 \d p  -A(q)\d p\\
\omega_{0,q}(p):= \le(\frac {\theta_1(p-q)}{\theta_1(p-q)}-\frac {\theta_1'(p)}{\theta_1(p)}\ri)
\eea
where the matrix $A(q)$ is determined uniquely by the condition that $M^{-1}(p)\Cauchy_0(p,q)$ is locally analytic for $p\in \scr T$. This gives the two linear systems:
\be
\label{Asys}
{\bf v}_j A(q) = \omega_{0,q}(t_j) {\bf v}_j
\ee
which can be solved  by 
\be
A(q) = P \le[\begin{array}{cc}
\omega_{0,q}(t_1) & 0 \\
0 & \omega_{0,q}(t_2)
\end{array}\ri] P^{-1}, \ \ \ \ \ P = \le[
\begin{array}{c}
{\bf v}_1\\
{\bf v}_2
\end{array}
\ri]\ .
\ee
Thus the existence rests uniquely upon the independence of the two vectors ${\bf v}_{1,2}$ established in Proposition \ref{proptyurin}.  

\subsubsection{Small norm theorem}
With the tool of the non-abelian Cauchy kernel at hand we can now state and formulate the existence and uniqueness of the solution of the RHP \ref{RRHP} for sufficiently large $n$.

In fact we can state and prove a more general result below.
In the formulation and proof we use the convention that for a matrix $X$ we denote by $|X|$ the Hilbert--Schmidt norm $(\tr (X \ov X^t))^\frac 1 2$, or any equivalent norm, and for a matrix valued function $X(p)$ we denote by  $\| X\|_s$ the $L^s$ norm of $|X(p)|$.
\bet
\label{smallR}
Let $\Gamma$ be a compact collection of smooth contours, locally finite on $\mathcal E_\tau$. Let $G(p)$ be a smooth matrix--valued function defined on $\Gamma$. Then the solution of the RHP
\bea
R(p_+) = R(p_-)\bigg(\1  +  G(p)\bigg),\ \ \ \ p\in \Gamma,\ \ \ R(0)=\1\\
\label{Rjump}
R(p) M(p) = \mathcal O(1), \ \  p\in \scr T
\eea
where $R$ is a matrix valued function analytic in $\E_\tau\setminus (\Gamma \cup \scr T)$, 
exists and is unique for $\|G\|_\infty$ sufficiently small. The matrix furthermore satisfies 
\be
| R(p)-\1| = \frac{\mathcal O(\|G\|_\infty)}{{\rm dist} (p, \Gamma \cup \scr T)}\label{Rest1} 
\ee
uniformly on compact sets in $\mathcal E_\tau \setminus (\Gamma \cup \scr T)$.  
More generally there is a constant $K$ such that 
\be
\label{Rest2}
| R(p)M(p)-M(p)| \leq \frac{ K\|G\|_\infty}{ {\rm dist}(p, \Gamma\cup\scr D)} ,
\ee
where the distance is in the metric $|\d p|$ (i.e. the flat metric on the torus). 

Assuming, furthermore that $G(p)$ is analytic in a tubular neighbourhood of $\Gamma$, then we have the stronger estimate 
\be
| R(p)-\1| \leq \frac{ K\|G\|_\infty}{{\rm dist} (p,\scr T)},\qquad 
| R(p)M(p)-M(p)| \leq \frac{ K\|G\|_\infty}{{\rm dist} (p,\scr D)},
\ee
for some constant $K$.
\eet
\noindent {\bf Proof.} 
We claim that the solution satisfies the integral equation:
\be
R(q)  = \1 + \frac 1{2i\pi} \int_\Gamma  R(p_-)  G(p) \Cauchy_{0}(p,q)\d p
\label{Volterra}
\ee
Indeed, since the singularity along the diagonal $p=q$ of the Cauchy kernel is of the form $\1/(p-q)$, it then  follows from the Sokhotski--Plemelj formula  that $R(p_+)-R(p_-)  =  R(p_-)G(p)$, which implies the jump condition \eqref{Rjump}. 
Moreover, since $\Cauchy_0(p,0)\equiv 0$, also the normalization $R(0)=\1$ is guaranteed. Finally the analyticity of $R(q)M(q)$ at the Tyurin points follows from the property that $\Cauchy_0(p,q)M(q)$ is analytic when $q\in \scr T$, by definition of non-abelian Cauchy kernel. 

The Cauchy integral operator allows us to take boundary values $q\in \Gamma$ from either sides and the result is a bounded operator in $L^2(\gamma, |\d p|)$.  
Taking the boundary value we obtain the integral equation 
\be
({\rm Id} - \mathcal G\big) [R_-](q)=  \1\label{inteq}
\ee
where  $\mathcal G$ acts on the each row ${\bf r}$ of $R_-$ as the following operator
\be
\mathcal G[{\bf r} ](q) = \int_{\gamma}{\bf r}(p) G(p) \Cauchy_0(p, q_-) \frac{\d p}{2i\pi}
\ee 
as an operator in $L^2(\Gamma, |\d p|) \otimes \C^2$. Let $N_1$ be the operator norm of the Cauchy boundary operator; then one finds immediately that  the norm of $\mathcal G$ is bounded by  $N_1 \|G\|_\infty$.
 Thus the operator ${\rm Id} -  \mathcal G$ is invertible if $\|G\|_\infty$ is smaller than $1/N_1$, and the inverse has also bounded operator norm  (denoted here by $\|\hspace{-0.8pt}|\bullet \|\hspace{-0.8pt}|$)
\be
\|\hspace{-0.8pt}|({\rm Id}-  \mathcal G)^{-1}\|\hspace{-0.8pt}|_2\leq \frac 1{1-\|\hspace{-0.8pt}|\mathcal G\|\hspace{-0.8pt}|_2} \leq\frac 1{1-N_1\|G\|_\infty}.
\ee
Thus the $L^2$ norm of $R(q_-)$ is a priori bounded  as long as $\|G\|_\infty$ is smaller than $N^{-1}_1$. 

Now consider $R(q)-\1$ for $q\not\in \Gamma$; from \eqref{Volterra} we have the point-wise identity
\be
\le |R(q)-\1\ri| = \le| \int_\Gamma R(p_-) G(p) \Cauchy_0(p,q)\frac{\d p}{2i\pi}  \ri|\ .
\ee
The non-abelian Cauchy kernel  $\Cauchy_0(p,q)$ is analytic for $p\in \Gamma$ if $q\not\in \Gamma$ and thus we have immediately 

\be
\big |R(q)-\1 \big | \leq \|R_-\|_{_2} \,\|G\|_{_2} \,\|\Cauchy_0(\bullet,q)\|_{_\infty}.
\ee
Here we mean that the $L^\infty$ norm of the matrix $\Cauchy_0$ is taken with respect to the first variable (and the result is thus a positive function of $q$).
Since $\Gamma$ is necessarily compact, the $L^2$ norm of $G$ is bounded above by the $L^\infty$ norm  times the length of $\Gamma$ (in the metric $|\d p|$) and hence 
\be
\le |R(q)-\1\ri | \leq \frac{ L(\Gamma) \|G\|_{_\infty}}{1-N_1\|G\|_\infty} \|\Cauchy_0(\bullet,q)\|_{_\infty}.
\ee 
Finally, since the singularity of $\Cauchy_0(p,q)$ along the diagonal $p=q$ is a simple pole with residue one, we can bound its matrix norm by a constant  over the distance from $\Gamma$ and its other poles with respect to $q$ at the Tyurin points $\scr T$ and we obtain \eqref{Rest1}.

To establish \eqref{Rest2} we proceed similarly 
\be
\big |R(q)M(q) -M(p)\big | \leq \|R_-\|_{_2} \,\|G\|_{_2} \,\|\Cauchy_0(\bullet,q)M(q) \|_{_\infty}.
\ee
Since now $\Cauchy_0(\bullet,q)M(q)$ is now bounded also for $q\in \scr T$ (according to the defining properties of the non-abelian Cauchy kernel), we can bound it by the reciprocal of the  distance from $\Gamma$ and the poles of $M(q)$ at $\scr D$  on $\mathcal E_\tau$ and we obtain \eqref{Rest2}.

To complete the proof we now assume furthermore the analyticity of $G(p)$ in a tubular neighbourhood of $\Gamma$. This implies that $R(p_-)$ is then the boundary value of an analytic function; thus we conclude that the estimate holds also on $\Gamma$ using the standard arguments.\QED

We apply the small norm Theorem \ref{smallR} to the RHP \ref{RRHP}; in this case 
\be
G(p;n)= {\rm e}^{-2\varphi(p)} M(p) \le[
\begin{array}{cc}
0 &0 \\ 1&0
\end{array}
\ri] M(p)^{-1}
\ee
and we have already established that 
\be
\sup_{p\in \gamma_\pm} \big |G(p;n)\big | \leq K {\rm e}^{-n c_0}
\ee
with $c_0>0$ given, for example, by \eqref{defc0}. 
This implies the following corollary.
\begin{corollary}
\label{Rissmall}
The solution $R(p)$ of the Problem \ref{RRHP} exists for any $n$. Moreover there is a constant $K>0 $ and $c_0>0$ (given by \eqref{defc0} for example) so that 
\be
\|R(p)M(p)- M(p)\| \leq \frac{K {\rm e}^{-nc_0}}{{\rm dist}( p, \scr D)}, \ \ \ \ 
\|R(p)- \1\| \leq \frac{K {\rm e}^{-nc_0}}{{\rm dist}( p, \scr T)}, \ \ \ \ 
\ee
\end{corollary}
{  
The only part where we may add a comment is the stated existence: since $R = T M^{-1}$ is defined in terms of the original Problem \ref{YRHP} and the model Problem   \ref{MRHP} (both of which have been shown to exist for any $n\in \mathbb N$) via a sequence of exact transformations,  then $R$ solves Problem \ref{RRHP} and hence the solution exists. 
}
\subsection{Epilogue: asymptotics and zero density}
\label{epilogue}
The orthogonal ``monic'' section $\pi_n(p)$ is the $(1,1)$ entry of the original matrix $Y$ and hence, tracing backwards the transformations, 
\bea
Y(p) &\mathop{=}^{\eqref{Zdef}}  {\rm e}^{((n-1) \ell - S_\infty) \s_3} Z {\rm e}^{((n-1) (g(p)-\ell) +S(p))\s_3} 
\nn 
\\
&\mathop{=}^{\eqref{Tdef}}   {\rm e}^{((n-1) \ell - S_\infty) \s_3} T
\le[\begin{array}{cc}
1 & 0 \\
\pm \chi_{\pm}{\rm e}^{-2\varphi} & 1
\end{array}\ri]
 {\rm e}^{((n-1) (g(p)-\ell)+S(p))\s_3} 
 \nn
 \\
& \mathop{=}^{\eqref{Rdef}}   {\rm e}^{((n-1) \ell - S_\infty) \s_3} R(p)M(p)
\le[\begin{array}{cc}
1 & 0 \\
\pm\chi_\pm  {\rm e}^{-2\varphi} & 1
\end{array}\ri]
 {\rm e}^{((n-1) (g(p)-\ell)+S(p))\s_3}
 \label{chain} 
\eea
where $\chi_{_{\pm}}$ are the characteristic functions of the regions $\mathcal L_\pm$ (see Fig. \ref{figlenses}). 
{
According to Corollary \ref{Rissmall} we can replace the term $R(p) M(p)$ with $M(p)$ up to addition of terms of order $\frac{\mathcal O( {\rm e}^{-nc_0})}{{\rm dist}( p, \scr D)}$. Since $\gamma \cap \scr D=\emptyset$ we can drop the denominator of this estimate in a small tubular neighbourhood of $\gamma$. 
}
\paragraph{Outside the strip-lens.}
For $p$ outside of the lens-strip  and recalling that $T(p)  = R(p) M(p)$ (with $T(p)$ given by \eqref{Tdef}, \eqref{Zdef}) we have the following asymptotics
{ 
\bea
\pi_n(p) =&   {\rm e}^{ - S_\infty}  \le( M_{11}(p)  +\frac{\mathcal O({\rm e}^{-n c_0})}{{\rm dist}( p, \scr D)} \ri)
{\rm e}^{(n-1) g(p) + S(p)}  \ .
\label{pionstrip}
\eea
}
If $n $ is even  then $M_{11}(p)$ has a zero at $p=\frac 12$ (outside of the lens-strip).
If $n $ is odd  then $M_{11}(p)$ has no zeros on $\E_\tau\setminus \gamma$. This implies that for even $n$ the extra zero of $\pi_n$ converges to $p=\frac  1 2 $ at exponential rate by the Rouch\'e\ theorem. 
Note the case where $\scr D  = (1/2)$; in this case $M_{11}$ is actually regular at $p=1/2$ and has no zeros;  thus we can still conclude that the zero of $\pi_n$ also tends to $(1/2)$ (and ``cancels'' the pole).

\paragraph{In the strip-lens (and on the contour).}
Now we can describe the  nature of the oscillatory behaviour of the sections on the contour $\gamma$. This is achieved by tracing again the chain of transformations. Consider the upper half of the strip-lens and the equality \eqref{chain}. Taking the $(1,1)$ entry in \eqref{chain} we obtain:
{ 
\bea
\pi_n 
&={\rm e}^{(n-1) \ell - S_\infty} \le(
M_{1 1}(p_+) {\rm e}^{(n-1)  (g(p_+)-\ell) + S(p_+) }  
+ 
M_{12}(p_+) {\rm e}^{(n-1)  (g(p_+) -\ell) +S(p_+)-2\varphi(p_+)}+ \mathcal O({\rm e}^{-nc_0})
\ri)\ .
\eea
}
Using  the definition of the effective potential $\varphi$ \eqref{defphi} we can rewrite the above as 
{ 
\bea
\pi_n  
&={\rm e}^{(n-1) \ell - S_\infty} {\rm e}^{-\frac 1 2 w(p) } \le(
M_{1 1}(p_+) {\rm e}^{\varphi(p_+)  }  
+ 
M_{12}(p_+) {\rm e}^{-\varphi(p_+)}+ \mathcal O({\rm e}^{-nc_0})
\ri) \ .
\eea
}
We use now that $\varphi(p_+)\in i\R$ and $\varphi(p_-)=-\varphi(p_+)$ together with the fact that $M_{12}(p_+) = \ov M_{11}(p_+)$ (from Corollary \ref{SchartzM}) to get 
{ 
\be
\pi_n(p) =2 {\rm e}^{(n-1) \ell - S_\infty} {\rm e}^{-\frac 1 2 w(p) }   \Re \bigg(M_{11}(p_+) {\rm e}^{\varphi(p_+)} + \mathcal O({\rm e}^{-nc_0})\bigg) \ .
\ee
}
Denote by $A(s) = \le|M_{11}\le(s+\frac \tau 2 +i0\ri)\ri|$ and $\rho(p)$ the phase: observe the since $M_{11}(p+1)=(-)^{n+1} M_{11}(p)$
the phase has an increment of $\pi $ modulo $2\pi$ for $(n-1) $ odd, thus $\rho$ is a function that depends on the parity of $n$. 
Then we can write 
{ 
\be
\pi_n(p) = 2 {\rm e}^{(n-1) \ell - S_\infty} {\rm e}^{-\frac 1 2 w(p) } \le[
A(p)  \cos \le(-(n-1) \int_{\frac \tau 2} ^p \sqrt{e_1-\wp(q)} \d q + \rho(p)  +  \nu(p) \ri) + \mathcal O({\rm e}^{-nc_0})\ri] \ ,
\ee}
where $\nu(p)$ is defined in \eqref{SSchwartz} and 
\be
\rho(p) = \arg(M_{11}(p)) = \arg\le(f(p)\ri)
\ee
with $f$ as in \eqref{defF}. This follows from the fact that the normalization $C_{11}$ \eqref{defC1}, \eqref{defC2} is real. Note that the argument of the cosine is a function with increment of multiples of $2\pi$ over the contour $\gamma$ both in the even and odd case.  

Using the explicit expression \eqref{explicitG} we derive 
\bea
\label{img}
 \Phi(s):= \Im g\le(s+ \frac \tau 2+  i0\ri)  = \pi s  -\frac \pi 2 +\arg \frac {\theta_1(s+\frac \tau 2;2\tau)}{\theta_1(s-\frac \tau 2;2\tau)}
\eea
and then we can write 
{  
\be
\pi_n\le(s+\frac \tau 2\ri) = 2 {\rm e}^{(n-1) \ell - S_\infty} {\rm e}^{-\frac 1 2 w(p) } \le[A(s)  \cos \le(
(n-1)\Phi(s)+ \rho\le(s+\frac \tau 2\ri)  +  \nu\le(s+\frac \tau 2\ri) \ri)+ \mathcal O({\rm e}^{-nc_0})\ri]\ .
\ee
}
The expression \eqref{approxintrosupp} follows from elementary simplifications.

\subsubsection{Density of zeros and convergence of the normalized counting measure}
Recall that all but possibly one of the zeroes of $\pi_n$ are on $\gamma$. 
\bet
The normalized density of zeroes of $\pi_n$  tends  weakly to the measure 
\be
\d\mu_0(p) =\frac 1 \pi  \sqrt{ \wp\le(\frac 1 2\ri) -\wp  (p) } \d p = \frac {\d X}{2\pi \sqrt{(e_2- X)(X-e_3)}}
\ee
whose density is shown in Fig. \ref{Figdensity}. The convergence is in the sense that for any  continuous function  $\phi$  we have 
\be
\lim_{n\to\infty}\frac 1 n \sum_{j=1}^{2\lfloor \frac {n+1}2 \rfloor} \phi(z_j^{(n)}) = \int_\gamma \phi(p) \d\mu_0(p).
\ee
where $\{z_1^{(n)},\dots, z_m^{(n)}\}$ are the zeros of $\pi_n$ on $\gamma$ (with $m =  2\lfloor \frac {n+1}2 \rfloor$).
\eet
{\bf Proof.}
We know  already that $\pi_n$ has all zeros on $\gamma$ except at most one  on $\alpha$ for even $n$.  
Suppose first  that $\phi$ is analytic in the strip $\le|\Im p-\frac 1 2 \Im \tau \ri|<\epsilon$ for some $\epsilon>0$. In particular $\phi(p+1)=\phi(p)$ since it is stated to be analytic in the strip within the elliptic curve $\mathcal E_\tau$. Then we can evaluate the function $\phi$ at the zeros of $\pi_n$ by means of the residue theorem with 
\be
\frac 1 n \sum_{j=1}^n \phi(z_j^{(n)})=
\frac 1{2i n \pi} \le[
\int_{\frac \tau 2 + i\epsilon}^{\frac \tau 2 + i\epsilon+1} \phi(p) \d \ln \pi_n(p)
-
\int_{\frac \tau 2 - i\epsilon+1}^{\frac \tau 2 - i\epsilon} \phi(p) \d \ln \pi_n(p)
\ri]\ .
\ee
We can then use \eqref{pionstrip} to conclude that 
\be
\lim_{n\to\infty} \frac 1 n \sum_{j=1}^n \phi(z_j^{(n)})=  
\frac 1{2i\pi} \le[
\int_{\frac \tau 2 + i\epsilon}^{\frac \tau 2 + i\epsilon+1} \phi(p) \d g(p)
-
\int_{\frac \tau 2 - i\epsilon+1}^{\frac \tau 2 - i\epsilon} \phi(p) \d g(p)
\ri].
\ee
We now can send $\epsilon\searrow 0$ without changing the value by the Cauchy theorem, and we obtain
\be
\lim_{n\to\infty} \frac 1 n \sum_{j=1}^n \phi(z_j^{(n)})=  
\frac 1{2i\pi} \int_{\gamma} \phi(p) \bigg(\d g(p_+ )-\d g(p_-)\bigg).
\ee
We then use the Definition \ref{defGfun} and also \eqref{3.8} to reach the conclusion. 

Now the proof for an arbitrary continuous function can be obtained using the usual density argument of analytic functions on $\gamma$ within the space of continuous functions (i.e. Stone-Weierstrass theorem). 
\QED
\subsubsection{Asymptotics of the norms}
\label{norms}
In order to compute the norms we observe that 
\be
\|\pi_n\|^2 =\int_{\gamma}\pi_n^2(p) \d \mu(p)=  \int_\gamma Y_{11}(p) \big(Y_{12}(p_+)- Y_{12}(p_-)\big) \d p\ . 
\ee
Now it follows from the formulation of the RHP \ref{YRHP} that the expression $Y_{11}(p) Y_{12}(p)$ is meromorphic on $\mathcal E_\tau \setminus \gamma$ with a   simple pole only at $p=0$ and we can thus write 
\be
\|\pi_n\|^2 =2i\pi\res{p=0} Y_{11}(p)Y_{12}(p) \d p\ .
\ee
This residue is simply a contour integral on a circle around the origin and hence we can use the asymptotics of $Y(p)$  from \eqref{chain}. A short computation using the formulas in Theorem \ref{thmM}  yields
\bea
2i\pi\res{p=0} Y_{11}(p) Y_{12}(p) &=2i\pi \res{p=0} {\rm e}^{2(n-1)\ell - 2S_\infty} M_{11}(p)M_{12}(p) \le(1 + \mathcal O({\rm e}^{-nc_0})\ri)
=\\
&=
2\pi  {\rm e}^{2(n-1)\ell - 2S_\infty} {\rm e}^{-i\pi \tau} 
\frac { {\rm e}^{-2i\pi \scr D} \theta_1^2(\scr D ;2\tau)}{\theta_1^2(\scr D+\tau;2\tau)}
\frac {\theta_1'(0;2\tau) }{\theta_4(0;2\tau) }
 \le(\frac { \theta_3(0;2\tau)}{\theta_2(0;2\tau)  }\ri)^{\sharp_n} \big(1 + \mathcal O({\rm e}^{-nc_0})\big)
\eea 
where $\sharp_n=1$ for $n$ even and $-1$ for $n$ odd.

\section{Conclusion and open problems}
\label{final}
The present paper constitutes the first example of nonlinear steepest descent analysis in higher genus, to the knowledge of the author. 
We first comment on some immediate generalization and then further directions.

The first question is  what happens if we waive  the requirement that $\mathcal E_\tau$ is a real elliptic curve. 
In this case
\begin{enumerate}
\item [-]In general the determinant $D_n$ is complex and could vanish for some $n$; hence for these values of $n$ the RHP \ref{YRHP} is not solvable or does not have a unique solution. 
\item [-]  Theorem \ref{thmzeros}  does not hold any longer and the location of the zeros cannot be pinned down. 
\item[-] The construction of the $g$ and Szeg\"o\ functions is identical; however, the steepest descent analysis hits an obstacle because the  vertical trajectories of the quadratic differential $(\d g)^2 = (\wp(p)-e_1) \d p^2$ on the elliptic curve in general do not close to a cycle. There are  a disk-domain where the trajectories close to circles homotopic to $p=0$,  and an open set where we have {\it recurrent} trajectories, namely, they fill densely an open subset of $\mathcal E_\tau$ \cite{Strebel}. This suggests that the zeros of the orthogonal section should fill densely an open region. However this is not entirely clear and further study is  needed. Unfortunately numerical results are inconclusive because the author was not able to produce sufficiently high-degree orthogonal section with the required numerical accuracy.
\end{enumerate}
While retaining the same general setup of real elliptic curves one may consider other situations where:
\begin{itemize}
\item[-] The support of the orthogonality measure ${\rm e}^{w(p)} \d p$ is not the whole $\gamma$. For example it is a single sub-arc $[a,b]$. This would affect the construction of the $g$--function which now has to satisfy $\d g_+ + \d g_- = $constant only on the support $[a,b]$ and hence $\d g$ will have inverse square-root singularities at the endpoints. This has the consequence that the differential $\d g$ lifts to a meromorphic differential on a double cover, $\wh {\mathcal C}$ of $\mathcal E_\tau$ branched along the support $[a,b]$ hence a Riemann surface of genus $g=2$. Similarly the solution of the corresponding model problem involves now Riemann Theta function associated with the genus $2$ cover $\wh {\mathcal C}$.
If the support has $k$ sub-arcs of $\gamma$, then this double cover will have genus $k+1$. 

{ 
 In this case one will also need to construct suitable local parametrices near the endpoints of the arcs; since the considerations are purely local, we expect the parametrices to be expressible as usual in terms of Bessel functions evaluated in terms of a suitable local coordinate. The use of local parametrices (be it of Bessel or Airy type) will have the immediate effect of reducing the rate of convergence of the approximation from exponential to algebraic decay in $n$. Of course the details deserve a separate study.} 
\item [-] The weight function depends on a scaling parameter $N$ as ${\rm e}^{Nw(p)} \d p$. In this case we would study the asymptotics of $\pi_n(p;N)$ as $n\to \infty$, $N\to \infty$ and the ratio tends to a positive constant. 
The general approach still applies but now the construction of the $g$ function  relies upon a version of the variational problem associated with an appropriate generalization to higher genus of the Green energy functional with external potential. 
\end{itemize}

More ambitiously one can venture in the land where the initial Riemann surface $\mathcal C$ is a (real) algebraic curve of genus $g\geq 2$. This requires some level of sophistication both in the construction of the $g$--function (be it in the fixed measure or scaling measure case) and in the construction of the corresponding model problem.  In this case the results of \cite{BNR} on the construction of the non-abelian Cauchy kernel will become even more fundamental and moreover the analysis is nicely intertwined with the theory of vector bundles. 

In a more speculative manner, we can imagine to apply the whole theory of orthogonality on Riemann surfaces (and ancillary problems of asymptotics)  to the same circle of problems where ordinary orthogonal polynomials are used, for example the construction of Determinantal Random Point Fields.

It would be also interesting to investigate connections to integrable systems if the weight ${\rm e}^{w(p)}$ depends on suitable deformation parameters, much in the same way as it happens in the realm of ordinary orthogonal polynomials and the KP and Toda hierarchies.

\appendix
\renewcommand{\theequation}{\Alph{section}.\arabic{equation}}
{  
\section{The scalar RHP in genus $1$}
\label{genex}
Let $\E_\tau$ be an elliptic curve and $\gamma$ a contractible loop, namely the boundary of an embedded disk $\DD$, e.g. the circle $|p|=\epsilon$ modulo $\Lambda_\tau$. Let $\gamma = \pa \DD$ and  $J:\gamma\to \C^\times$ be a smooth function of index zero, namely, such that $\ln J$ is defined and continuous on $\gamma$.
\bp
\label{propRHPmult}
Let the above assumptions  prevail. 
Consider the  (scalar) RHP
\begin{equation}
\label{multRHP}
Y(p_+)  = Y(p_-) J(p)\ ,\ \ p\in \gamma,
\end{equation}
with $Y(p)$ and $Y(p)^{-1}$ both  analytic in $\E_\tau\setminus \gamma$ and bounded everywhere. Then this problem 
has solution if and only if 
\begin{equation}
\label{condRHPmult}
\frac 1{2i\pi}\int_{\gamma} \ln J(p) \d p =m+\tau n , \ \ n,m\in \Z.
\end{equation}
\ep
\noindent{\bf Proof.}
Suppose $Y(p)$ is a solution of the problem and denote $Y_+$ the restriction to the inside of the disk with boundary $\gamma$ and $Y_-$ the restriction to the outside of the disk. Neither of the two functions has any zeros or poles and hence the logarithm is a well defined function on the universal cover of $\mathcal C$: in particular $\ln Y_-(p)$ may  gain an additive integer multiple of $2i\pi$ as we analytically continue it on $\E_\tau \setminus \DD$.
Since the index of $J$ is zero, $j(z):=  \ln J(z)$ is a well defined continuous function $j:\gamma \to \C$. 
Letting $y_\pm =  \ln Y_\pm$, we need to solve an additive version $y_+(p) - y_-(p)= j(p)$ with $p\in \pa \DD$.
Now multiply by $\d p$ both sides of this equation and integrate on the boundary $\pa \DD$. Since $y_+$ is analytic in $\DD$ the integral gives zero contribution; now, by Cauchy's theorem, 
\be
\int_\gamma y_-(p)\d p = \int_{\pa \mathcal L} y_-(p)\d p
\ee
where $\mathcal L$ is any fundamental domain (parallelogram) of $\C/(\Z+\tau \Z)$. 
The function $y_-$ on two opposite sides of such parallelogram must differ by an integer multiple of $2i\pi$. Thus we obtain 
\be
-\frac 1{2i\pi}\int_{\gamma} \ln J(p) \d p =\frac 1{2i\pi}\int_{\gamma} y_-(p)\d p =\frac 1{2i\pi}\int_{\pa \mathcal L} y_-(p)\d p = m + \tau n.
\ee
This concludes the proof of the necessity of the condition \eqref{condRHPmult}.

To prove the sufficiency, let $p_0\in \E_\tau \setminus \gamma$ be an arbitrary point and consider the integral 
\begin{equation}
\label{defu}
u(p) =\frac 1{2i\pi} \int_{q\in\gamma}  \ \le(\frac{\d}{\d q}  \ln \frac{\theta_1(q-p)}{\theta_1(q-p_0)}\ri) \ln J(q)\d q
\end{equation}
where $\theta_1$ is the Jacobi theta function. The expression 
$$\omega_{p,p_0}(q):= \le(\frac{\d}{\d q}  \ln \frac{\theta_1(q-p)}{\theta_1(q-p_0)}\ri)\d q
$$
 is the normalized Abelian differential of the third kind; it has a pole at $q=p$ with residue $+1$ and a pole at $q=p_0$ with residue $-1$ and the integral along the $a$--cycle (i.e. $q\to q+1$) vanishes. Using the periodicity properties of $\theta_1$ one verifies that 
 \be
 \omega_{p+1,p_0}(q) = \omega_{p,p_0}(q), \ \ \ \omega_{p+\tau,p_0}(q) = \omega_{p,p_0}(q) + 2i\pi
 \ee
Using the Sokhotski--Plemelj formula we see that the  scalar function $u$ defined in \eqref{defu} satisfies $u(p_+)-u(p_-) = \ln J(p)$, $p\in \pa \DD$. 
It has also the following monodromy 
\begin{equation}
u(p+1) = u(p),\ \ \ u(p+\tau) = u(p) +  \int_{q\in\gamma} \ln J(q){\d q}{}.
\end{equation}
Now suppose 
\begin{equation}
\frac 1{2i\pi}\int_{\gamma} \ln J(q) \d q=  m +  \tau  n,\ \  \ \ \  m,  n \in \Z.
\end{equation}
Then we propose the solution 
\begin{equation}
y(p) = u(p) -  2i\pi n \,p.
\end{equation}
Indeed now we have 
\begin{equation}
y(p+1) = -2i\pi  n\ ,\ \ \ y(p+\tau) = y(p) + \int_{\gamma} \ln J(q)\frac{\d q}{2i\pi} - 2i\pi \tau n = y(p) +  2i\pi m.
\end{equation}
So indeed $Y(p) = {\rm e}^{y(p)} $  defines by restriction to $\DD$ and $\E_\tau \setminus \DD$ two analytic nonzero functions $Y_\pm$, and in particular $Y_-$ is single--valued on the non-simply connected domain  $\E_\tau \setminus \DD$.
Note that the proposed solution satisfies additionally $Y(p_0)=1$. 
\QED

}

\end{document}